\documentclass{amsart}
\usepackage{amsthm}
\usepackage{amssymb}
\usepackage{upref}
\usepackage{amscd}
\usepackage{xspace}

\numberwithin{equation}{subsection}

\newtheorem{theorem}[equation]{Theorem}
\newtheorem{lemma}[equation]{Lemma}
\newtheorem{proposition}[equation]{Proposition}
\newtheorem{corollary}[equation]{Corollary}

\newtheorem{question}[equation]{Question}
\newtheorem{lemmasec}{Lemma}

\numberwithin{lemmasec}{section}
\newtheorem{theoremsec}[lemmasec]{Theorem}
\newtheorem{propositionsec}[lemmasec]{Proposition}
\newtheorem{corollarysec}[lemmasec]{Corollary}

\theoremstyle{definition}
\newtheorem{definition}[equation]{Definition}
\newtheorem*{remark}{Remark}

\newtheorem{definitionsec}[lemmasec]{Definition}

\DeclareMathOperator{\Hom}{Hom}
\DeclareMathOperator{\sho}{ho}
\DeclareMathOperator{\Ext}{Ext}

\DeclareMathOperator{\map}{map}
\DeclareMathOperator{\colim}{colim}

\DeclareMathOperator{\dom}{dom}
\DeclareMathOperator{\codom}{codom}

\DeclareMathOperator{\im}{im}

\DeclareMathOperator{\coker}{coker}

\renewcommand{\smash}{\wedge}
\newcommand{\boxprod}{\mathbin\square}

\newcommand{\SSet}{\mathbf{SSet}}

\newcommand{\Ch}[1]{\text{Ch} (#1)}
\newcommand{\stable}[1]{\text{Stable} (#1)}

\newcommand{\hopfalg}{(A,\Gamma )}
\newcommand{\otherhopfalg}{(B,\Sigma )}

\newcommand{\cat}[1]{\mathcal{#1}}

\newcommand{\Z}{\mathbb{Z}}
\newcommand{\Fp}{\mathbb{F}_{p}}
\newcommand{\Zp}{\Z _{(p)}}
\newcommand{\comod}{\text{-comod}}
\newcommand{\Mod}{\text{-mod}}

\newcommand{\ev}{\textup{Ev}}

\newcommand{\llp}{left lifting property with respect to\xspace}
\newcommand{\rlp}{right lifting property with respect to\xspace}

\newcommand{\mathcolon}{\colon\,}

\newcommand{\uc}{\textup{:}}

\newcommand{\ulp}{\textup{(}}
\newcommand{\urp}{\textup{)}}

\begin{document}
 
\title{Homotopy theory of comodules over a Hopf algebroid} 

\date{\today}

\author{Mark Hovey}
\address{Department of Mathematics \\ Wesleyan University
\\ Middletown, CT 06459}
\email{hovey@member.ams.org}


\begin{abstract}
Given a good homology theory $E$ and a topological space $X$, $E_{*}X$
is not just an $E_{*}$-module but also a comodule over the Hopf
algebroid $(E_{*}, E_{*}E)$.  We establish a framework for studying the
homological algebra of comodules over a well-behaved Hopf algebroid
$\hopfalg $.  That is, we construct the derived category $\stable{\Gamma
}$ of $\hopfalg $ as the homotopy category of a Quillen model structure
on $\Ch{\Gamma }$, the category of unbounded chain complexes of $\Gamma
$-comodules.  This derived category is obtained by inverting the
homotopy isomorphisms, \textbf{NOT} the homology isomorphisms.  We
establish the basic properties of $\stable{\Gamma }$, showing that it is
a compactly generated tensor triangulated category.  
\end{abstract}

\maketitle

\section*{Introduction}

Given a commutative ring $k$, a \emph{Hopf algebroid} over $k$ is a
cogroupoid object in the category of commutative $k$-algebras.  That is,
a Hopf algebroid is a pair $\hopfalg$ of commutative $k$-algebras, so
that, given a commutative $k$-algebra $R$, the set $k\text{-alg}(A,R)$
is naturally the objects of a groupoid with morphisms
$k\text{-alg}(\Gamma ,R)$.  This gives several structure maps of which
we remined the reader below.  The reason for our interest in Hopf
algebroids is that, if $E_{*}(-)$ is a well-behaved homology theory on
topological spaces, then $E_{*}X$ is naturally a comodule over the Hopf
algebroid $(E_{*}, E_{*}E)$, though Hopf algebroids have also arisen in
algebraic geometry in connection with stacks~\cite{faltings-chai}.  In
particular, the study of comodules over the Hopf algebroid $BP_{*}BP$
led to the Landweber exact functor theorem~\cite{land-exact}, a result
of fundamental importance in algebraic topology.  When $(E_{*}, E_{*}E)$
is a Hopf algebroid, the $E_{2}$-term of the Adams spectral sequence
based on $E$ is the bigraded $\Ext $ in the category of
$E_{*}E$-comodules.

Thus we would like to understand the homological algebra of comodules
over a Hopf algebroid.  The simplest kind of Hopf algebroid is a
\emph{discrete} Hopf algebroid $(A,A)$.  The associated groupoid has no
non-identity maps, and a comodule over $(A,A)$ is the same thing as an
$A$-module.  One of the most useful tools in studying the homological
algebra of $A$-modules is the unbounded derived category $\cat{D}(A)$,
obtained by inverting the homology isomorphisms in the category $\Ch{A}$
of unbounded chain complexes of $A$-modules.  The goal of this paper is
to construct $\cat{D}\hopfalg $, the derived category of a Hopf
algebroid $\hopfalg $.  We stress that homology isomorphisms are
\textbf{NOT} the right thing to invert to form $\cat{D}\hopfalg $.  This
is already clear in case $\hopfalg $ is a Hopf algebra over a field $k$,
such as the Steenrod algebra.  In this case, $\cat{D}\hopfalg $ was
constructed by the author in~\cite[Section~2.5]{hovey-model} and studied
by Palmieri~\cite{palmieri-book}.  The idea is that chain complexes of
comodules are like topological spaces; they have homotopy as well as
homology, and it is the homotopy isomorphisms we should invert, not the
homology isomorphisms.  For a discrete Hopf algebroid $(A,A)$, homotopy
and homology coincide, but not in general.  To avoid confusion, we refer
to $\cat{D}\hopfalg $ as the \emph{stable homotopy category of $\Gamma
$-comodules}, and denote it by $\stable{\Gamma }$.

We want $\stable{\Gamma }$ to have all the usual properties of
$\cat{D}(A)$ or the stable homotopy category; it should be a
triangulated category with a compatible closed symmetric monoidal
structure, and there should be a good set of generators.  In fact, we
want $\stable{\Gamma }$ to be a stable homotopy category in the sense
of~\cite{hovey-axiomatic}.  We also want 
\[
\stable{\Gamma }(S^{0}A,S^{0}M)_{*} \cong \Ext ^{*}_{\Gamma }(A,M)
\]
for a comodule $M$, where $S^{0}N$ denotes the complex consisting of $N$
in degree $0$ and $0$ everywhere else.  This will guarantee that we
recover the $E_{2}$-term of the Adams spectral sequence based on $E$.

The axioms for a stable homotopy category, or, indeed, even for a
triangulated category, are so painful to check that the best way to
construct such a category is as the homotopy category of a Quillen model
structure~\cite{quillen-htpy}.  One of the main goals
of~\cite{hovey-model} was to enumerate the conditions we need on a model
structure so that its homotopy category is a stable homotopy category.  
We also point out that there are many advantages of a model structure
over its associated homotopy category; the model structure allows one to
perform constructions, such as homotopy limits and colimits, that are
inaccessible in the homotopy category, and allows one to make
comparisons with other model categories.  

Thus, the bulk of this paper is devoted to constructing a model
structure on $\Ch{\Gamma }$ in which the weak equivalences are the
homotopy isomorphisms.  We expect that the associated stable homotopy
category $\stable{\Gamma }$ will have many good properties and will
provide insight into homotopy theory, as Palmieri's
work~\cite{palmieri-book,palmieri-quillen-strat} on the Steenrod algebra
has done.  We have in fact shown that $\stable{E_{*}E}$ is a Bousfield
localization of $\stable{BP_{*}BP}$ for any Landweber exact commutative
ring spectrum $E$ in~\cite{hovey-barcelona}; this then gives rise to a
general change of rings theorem containing the change of rings theorem
of Miller-Ravenel~\cite{miller-ravenel} and the author and
Sadofsky~\cite{hovey-sadofsky-picard}.  The construction of the model
structure is complicated enough that we do not discuss such applications
in this foundational paper.  We do establish some beginning properties
of $\stable{\Gamma }$ in Section~\ref{sec-stable}.  In particular, we
show that $\stable{\Gamma }$ is monogenic when $\Gamma =BP_{*}BP$ or
$\Gamma =E_{*}E$ for $E$ any Landweber exact homology theory over $BP$.

In order to establish our model structure, we need to first 
study the structure of the abelian category $\Gamma \comod $
of $\Gamma $-comodules, which we do in Section~\ref{sec-comodules}.
Most of the results in this section seem to be new, at least in the
generality in which we give them, and of independent interest.  For
example, we study duality in $\Gamma \comod $, showing that a comodule
is dualizable if and only if it is finitely generated and projective as
an $A$-module.  

We follow the usual plan to construct our model structure.  That is, we
start by building an auxiliary model structure in
Section~\ref{sec-model-proj} called the projective model structure.
This model structure is easy to construct, but has too few weak
equivalences (unless the Hopf algebroid is discrete).  So we must
localize it by making the homotopy isomorphisms weak equivalences.  This
first necessitates a study of the homotopy isomorphisms in
Section~\ref{sec-groups} and a reminder, with a few new results, about
localization of model categories in~\ref{sec-model}.  We finally
construct the desired model structure in~\ref{sec-homotopy}, and study
some of the basic properties of $\stable{\Gamma }$ in the aforementioned
Section~\ref{sec-stable}.  

We should note that our results do not apply to an arbitrary Hopf
algebroid.  We need our Hopf algebroid to be \emph{amenable}, defined
precisely in Definition~\ref{defn-amenable}.  All of the amenable Hopf
algebroids we know are in fact Adams Hopf algebroids, defined
in~\cite{goerss-hopkins-comodules} but implicit in Adams' blue
book~\cite[Section~III.13]{adams-blue}.  If $E$ is a ring spectrum that
satisfies Adams' condition, which we call \emph{topologically flat},
that $E$ be a minimal weak colimit of finite spectra $X_{\alpha }$ such
that $E_{*}X_{\alpha }$ is finitely generated and projective over
$E_{*}$, then $(E_{*}, E_{*}E)$ is an Adams Hopf algebroid
(Section~\ref{subsec-Adams}).  To make sure our results apply in cases
of interest, we must check that interesting ring spectra $E$ are
topologically flat.  We prove in Theorem~\ref{thm-Adams-Landweber} that
if $R$ is topologically flat and $E$ is Landweber exact over $R$, then
$E$ is topologically flat.  

The author has been trying to prove the results in this paper since
1997, when Doug Ravenel strongly encouraged him to build a stable
homotopy category of $BP_{*}BP$-comodules.  It is a pleasure to
acknowledge the author's debt to Neil Strickland, who constructed
$\stable{BP_{*}BP}$ in a fairly ad hoc way, without a model structure,
about 1997.  The crucial input that finally enabled the author to build
the model structure came from the paper of Paul Goerss and Mike
Hopkins~\cite{goerss-hopkins-comodules}.

\section{The abelian category of comodules}\label{sec-comodules}

We begin with a fairly comprehensive study of the category $\Gamma
\comod $ of comodules over a Hopf algebroid $\hopfalg $.  Some of these
results are well-known, but others are apparently new.  

Before we begin, we establish notation and remind the reader of some of
the basic structure maps of Hopf algebroids.  The symbol $\hopfalg $
will always denote a Hopf algebroid~\cite[Appendix~1]{ravenel}, and the
symbol $\otimes $ always denotes $\otimes _{A}$, the tensor product of
$A$-\textbf{bimodules}.  Given an $A$-bimodule $M$, $\widetilde{M}$
denotes $M$ with the $A$-actions reversed.  

With these conventions, the structure maps of $\hopfalg $ include maps
of commutative $k$-algebras $\eta _{L}\mathcolon A\xrightarrow{}\Gamma $
corepresenting the source of a morphism, $\eta _{R}\mathcolon
A\xrightarrow{}\Gamma $ corepresenting the target of a morphism, and
$\varepsilon \mathcolon \Gamma \xrightarrow{}A$ corepresenting the
identity maps of the groupoid.  This makes $\Gamma $ into an
$A$-\textbf{bimodule}, with $\eta _{L}$ giving the left $A$-action and
$\eta _{R}$ giving the right $A$-action.  There are then additional
structure maps of $k$-algebras $\chi \mathcolon \Gamma
\xrightarrow{}\widetilde{\Gamma }$ corepresenting the inverse of a
morphism, and $\Delta \mathcolon \Gamma \xrightarrow{} \Gamma \otimes
\Gamma $ corepresenting the composition of a pair of morphisms.  Of
course, these maps must satisfy some relations assuring that we get a
groupoid.  For example, $\varepsilon \eta _{R}=\varepsilon \eta
_{L}=1_{A}$, since the source and target of the identity map at $x$ are
both $x$.  The remaining relations can be found
in~\cite[Appendix~1]{ravenel}.

\subsection{Basic structure}\label{subsec-comod-basic}

Recall that a \emph{left $\Gamma $-comodule} is a left $A$-module $M$
equipped with a map $\psi \mathcolon M \xrightarrow{}\Gamma \otimes M$
satisfying a coassociativity and counit condition.  There is an obvious
notion of a map of comodules.  

\begin{lemma}\label{lem-abelian}
Suppose $\Gamma $ is flat as a right $A$-module.  Then the category
$\Gamma \comod $ is a cocomplete abelian subcategory of $A\Mod$.  
\end{lemma}

\begin{proof}
Since the tensor product commutes with colimits, the $A$-module colimit
of a diagram of comodules is again a comodule, and is the colimit in
$\Gamma \comod $.  That $\Gamma \comod $ is abelian when $\Gamma $ is
flat is proved in~\cite[Theorem~A1.1.3]{ravenel}; we require flatness in
order to conclude that the $A$-module kernel of a comodule map is again
a comodule.
\end{proof}

Because of this lemma, we will assume throughout the paper that
$\hopfalg $ is a \emph{flat} Hopf algebroid; that is, that $\Gamma
$ is flat as a right $A$-module.  Note that the conjugation $\chi $
defines an isomorphism between the left $A$-module $\Gamma $ and the
right $A$-module $\Gamma $, so $\Gamma $ is also flat as a left
$A$-module.  

\begin{lemma}\label{lem-comod-monoidal}
The category $\Gamma \comod $ is a symmetric monoidal category. We
denote the symmetric monoidal product by $M\smash N$.  
\end{lemma}

\begin{proof}
We define $M\smash N=M\otimes N$, the tensor product of \textbf{left}
$A$-modules, with comodule structure given by the composite
\[
M\otimes N \xrightarrow{\psi \otimes \psi } (\Gamma \otimes M) \otimes
(\Gamma \otimes N) \xrightarrow{g} \Gamma \otimes M\otimes N,
\]
where $g(x\otimes m\otimes y\otimes n)=xy \otimes m \otimes n$.  We
leave it to the reader to check that this does define a map from the
tensor product, and that the composition above is a comodule structure.
Note that the multiplication map $\mu \mathcolon \Gamma \otimes
_{k}\Gamma \xrightarrow{}\Gamma $ does not factor through $\Gamma
\otimes \Gamma $, and this is why we use $g$.  The unit of the tensor
product is $A$, with comodule structure given by $\eta _{L}$.  
\end{proof}

We now point out that the category of comodules is natural.  
Recall that a map $\Phi \mathcolon \hopfalg \xrightarrow{}
\otherhopfalg$ of Hopf algebroids is a pair of ring homomorphisms
$\Phi_{0}\mathcolon A\xrightarrow{}B$ and $\Phi_{1}\mathcolon \Gamma
\xrightarrow{}\Sigma$ that corepresents a natural morphism of
groupoids.  This means that $\Phi_{0}\epsilon =\epsilon \Phi_{1}$,
$\Phi_{1}\eta_{L}=\eta_{L}\Phi_{0}$,
$\Phi_{1}\eta_{R}=\eta_{R}\Phi_{0}$, and $(\Phi_{1}\otimes
\Phi_{1})\Delta =\Delta \Phi_{1}$.  

\begin{lemma}\label{lem-naturality}
A map $\Phi \mathcolon \hopfalg \xrightarrow{}\otherhopfalg $ induces a
symmetric monoidal functor $\Phi _{*}\mathcolon \Gamma \comod
\xrightarrow{}\Sigma \comod $.
\end{lemma}

\begin{proof}
Define $\Phi_{*}M =B\otimes_{A}M$.  The $\Sigma $-comodule structure on
$B\otimes_{A}M$ is given by the composite
\[
B\otimes M \xrightarrow{1\otimes \psi} B\otimes \Gamma \otimes M
\xrightarrow{} \Sigma \otimes M\cong \Sigma \otimes_{B}(B\otimes M),
\]
where the map $B\otimes \Gamma \xrightarrow{}\Sigma$ takes $b\otimes x$
to $b\Phi_{1}(x)$.  
\end{proof}

In light of this lemma, the following definition is natural. 

\begin{definition}\label{defn-weak-equiv}
A map $\Phi \mathcolon \hopfalg \xrightarrow{}\otherhopfalg $ is a
\emph{weak equivalence} if $\Phi _{*}$ is an equivalence of categories.  
\end{definition}

Note that weak equivalences are a fundamentally new feature that arises
in studying Hopf algebroids; any weak equivalence of discrete Hopf
algebroids is necessarily an isomorphism, but there are many examples of
weak equivalences of Hopf algebroids that are not isomorphisms given
in~\cite{hovey-hopf} and~\cite{hovey-strickland-comodules}.  The author
used a different definition of weak equivalence in~\cite{hovey-hopf},
but the two definitions are in fact
equivalent~\cite{hovey-strickland-comodules}.

A particular example of a map of Hopf algebroids is the map $\epsilon
\mathcolon \hopfalg \xrightarrow{}(A,A)$ that is the identity on $A$ and
the counit $\epsilon $ on $\Gamma $.  Geometrically, this is the
inclusion of the identity maps of a groupoid into the whole groupoid.
The functor $\epsilon _{*}$ is just the forgetful functor from $\Gamma
$-comodules to $A$-modules.  As is well known~\cite[A1.2.1]{ravenel},
this functor has a right adjoint that takes an $A$-module $M$ to the
$\Gamma$-comodule $\Gamma \otimes M$, with structure map $\Delta \otimes
1$.  This is called the \emph{extended comodule} on $M$; in case $M$ is
itself a free $A$-module on the set $S$, then $\Gamma \otimes M$ is
called the \emph{cofree} comodule on $M$.  We have a natural isomorphism
\[
A\Mod (M,N)\xrightarrow{}\Gamma \comod (M,\Gamma \otimes N)
\]
for $\Gamma$-comodules $M$ and $A$-modules $N$.  This natural
isomorphism takes a map $f\mathcolon M\xrightarrow{}N$ of $A$-modules to
the map of comodules $(1\otimes f)\psi$, and a map of comodules
$g\mathcolon M\xrightarrow{} \Gamma \otimes N $ to the map $(\epsilon
\otimes 1)g$ of $A$-modules.  

It is less well-known that the extended comodule functor $M\mapsto
\Gamma \otimes M$ itself has a right adjoint $R\mathcolon \Gamma \comod
\xrightarrow{} A\Mod$, defined by $RN=\Gamma \comod (\Gamma ,N)$.  The
$A$-module action on $RN$ is defined by $(af)(x)=f(x\eta_{R}(a))$.  

Note that, if $M$ is itself a comodule, then we can form the extended
comodule $\Gamma \otimes M$ and the tensor product $\Gamma \smash M$.
The following lemma is well-known. 

\begin{lemma}\label{lem-different-tensor}
Suppose $\hopfalg$ is a flat Hopf algebroid, $M$ is an $A$-module, and
$N$ is a $\Gamma$-comodule.  Then there is a natural isomorphism of
comodules 
\[
(\Gamma \otimes M) \smash N \xrightarrow{} \Gamma \otimes (M\otimes N).  
\]
In particular, when $M=A$, we get a natural isomorphism of comodules 
\[
\Gamma \smash N\xrightarrow{} \Gamma \otimes N.
\]
\end{lemma}

\begin{proof}
We first note that $(\Gamma \otimes M)\smash N$ is the tensor product of
the \textbf{left} $A$-modules $\Gamma \otimes M$ and $N$.  There is a
natural comodule map 
\[
f_{MN}\mathcolon (\Gamma \otimes M)\smash N\xrightarrow{}\Gamma
\otimes (M\otimes N)
\]
adjoint to $\varepsilon \otimes 1\otimes 1$.  For fixed $N$, this is a
natural transformation of right exact functors of $M$ that commutes with
direct sums.  Since every $A$-module is a quotient of a map of free
$A$-modules, it suffices to show that $f_{AN}$ is an isomorphism.  

In fact, we construct an inverse $g=g_{AN}$ to $f_{AN}$.  
We define $g$ to be the composite 
\[
\Gamma \otimes M \xrightarrow{1\otimes \psi } \Gamma \otimes \Gamma
\otimes M \xrightarrow{(\mu\circ (1\otimes \chi ))\otimes 1} \Gamma \smash 
M,
\]
which is, \emph{a priori}, only a map of $A$-modules.  Note that, though
the multiplication $\mu $ does not factor through $\Gamma \otimes \Gamma
$, the composite $\mu \circ (1\otimes \chi )$ does do so, since $\chi $
switches the left and right units.  A diagram chase shows that $g$ and
$f_{AN}$ are inverses (and therefore that $g$ is a comodule map).
\end{proof}

It is tempting to think that, given an $A$-module $M$, one can think of
$M$ as a trivial $\Gamma$-comodule, via the map $\eta_{L}\otimes
1\mathcolon M \xrightarrow{} \Gamma \otimes M$.  This is wrong; for
example, $v_{n}^{-1}BP_{*}$ cannot be given the structure of a
$BP_{*}BP$-comodule~\cite[Proposition~2.9]{johnson-yosimura}.  The
difficulty is that $\eta_{L}$ is not a map of $A$-bimodules.  However,
there is a symmetric monoidal trivial comodule functor from the category
of abelian groups to $\Gamma$-comodules that takes the abelian group $M$
to $A\otimes_{\Z}M$ with the trivial comodule structure given by
$\eta_{L}\otimes_{\Z} 1$.  This functor has a right adjoint that takes
the comodule $N$ to the abelian group of primitive elements in $N$.

\subsection{Limits}\label{subsec-comod-limits}

In general, right adjoints such as limits are difficult to construct for
comodules, because the forgetful functor from $\Gamma $-comodules to
$A$-modules does not preserve products, though it does preserve kernels.
We give a general (and apparently new) method for constructing right
adjoints, involving resolutions by extended comodules.

For a comodule $M$, the adjoint to the identity map is the map $\psi
\mathcolon M\xrightarrow{}\Gamma \otimes M$, which we now think of as a
map of comodules, giving $\Gamma \otimes M$ the extended comodule
structure.  The map $\psi$ is of course an embedding, since it is split
over $A$ by $\epsilon \otimes 1$.  In particular, if $p\mathcolon \Gamma
\otimes M\xrightarrow{}N$ denotes the cokernel of $\psi$, which is
itself a comodule, then we have a natural diagram
\begin{equation}\label{eq-comod}
M\xrightarrow{\psi} \Gamma \otimes M \xrightarrow{\psi p} \Gamma \otimes N 
\end{equation}
expressing $M$ as the kernel of a map of extended comodules.  

Now, if $R$ is a right adjoint, then $R$ will have to preserve kernels,
so $R$ is completely determined by its restriction to the full
subcategory of extended comodules.  

We first use this idea to show that $\Gamma \comod$ is complete.  

\begin{proposition}\label{prop-comod-products}
Suppose $\hopfalg$ is a flat Hopf algebroid.  Then $\Gamma \comod$ has
products, and so is complete.  
\end{proposition}

\begin{proof}
Let us denote the comodule product we are trying to construct by
$\prod_{i}^{\Gamma} M_{i}$.  Adjointness shows that, if $\{M_{i} \}$ is
a set of $A$-modules, then 
\[
\prod_{i}^{\Gamma }(\Gamma \otimes M_{i})\cong \Gamma \otimes
\prod_{i}M_{i}.
\]
Now suppose we have a set of comodule maps $f_{i}\mathcolon \Gamma
\otimes M_{i}\xrightarrow{}\Gamma \otimes N_{i}$.  We need to define the
comodule product $\prod_{i}^{\Gamma}f_{i}$.  We define it to
be the composite
\begin{gather*}
\Gamma \otimes \prod M_{i} \xrightarrow{\Delta \otimes 1} \Gamma \otimes
\Gamma \otimes \prod M_{i} \xrightarrow{1\otimes \alpha} \Gamma \otimes
\prod (\Gamma \otimes M_{i}) \\
\xrightarrow{1\otimes \prod f_{i}} \Gamma \otimes \prod (\Gamma \otimes
N_{i}) \xrightarrow{1\otimes \prod (\epsilon \otimes 1)} \Gamma \otimes
\prod N_{i},
\end{gather*}
where $\alpha$ is the evident natural transformation.  The reader can
then check that this is a good definition of the product on the full
subcategory of extended comodules, and, in particular, that
$\prod_{i}^{\Gamma}(1\otimes g_{i})=1\otimes \prod g_{i}$ for $A$-module
maps $g_{i}$.

The definition of the product for a family of general comodules $M_{i}$
is now forced on us, as explained in the paragraph preceding this
proposition.  To wit, given a set of comodules $M_{i}$, we have left
exact sequences 
\[
0\xrightarrow{} M_{i} \xrightarrow{\psi} \Gamma \otimes M_{i}
\xrightarrow{\psi p_{i}} \Gamma \otimes N_{i}, 
\]
and so we define $\prod_{i}^{\Gamma}M_{i}$ by the left exact sequence 
\[
0 \xrightarrow{} \prod_{i}^{\Gamma} M_{i} \xrightarrow{} \Gamma \otimes
\prod_{i} M_{i} \xrightarrow{\prod_{i}^{\Gamma}(\psi p_{i})} \Gamma
\otimes \prod_{i} N_{i}.  
\]
We leave to the reader the proof that $\prod_{i}^{\Gamma}M_{i}$ is
indeed the product in $\Gamma \comod$.  
\end{proof}

\begin{remark}
An alternative approach to the category of $\Gamma$-comodules that is
sometimes used (e.g., by Boardman~\cite{boardman-eight}) is to establish
an equivalence of categories between $\Gamma$-comodules and a
subcategory of $\Gamma^{*}$-modules.  Here
$\Gamma^{*}=\Hom^{r}_{A}(\Gamma ,A)$, the $A$-bimodule of right
$A$-module maps from $\Gamma$ to $A$.  It turns out that $\Gamma^{*}$ is
a (noncommutative) algebra over $k$ and there is a map of algebras
$A\xrightarrow{}\Gamma^{*}$.  There is a map $\Gamma \otimes
M\xrightarrow{\alpha }\Hom_{A}(\Gamma^{*},M)$.  Using this map, any
$\Gamma$-comodule becomes a $\Gamma^{*}$-module, and this clearly
defines a faithful functor from $\Gamma$-comodules to
$\Gamma^{*}$-modules.  However, this functor will not in general be
full, because the map $\alpha$ need not be injective.  If $\Gamma$ is
projective over $A$, then $\alpha$ is injective, but $\Gamma $ is not
projective over $A$ for many of the Hopf algebroids we are interested
in.  If $\Gamma $ is projective over $A$, one can establish an
equivalence between $\Gamma$-comodules and a full coreflective
subcategory of $\Gamma^{*}$-modules.  This means that the inclusion
functor from $\Gamma $-comodules to $\Gamma ^{*}$-modules has a right
adjoint $R$.  Indeed, if $M$ is a $\Gamma ^{*}$-module, let us denote by
$\mu ^{*}\mathcolon M\xrightarrow{}\Hom _{A}(\Gamma ^{*},M)$ the adjoint
to the structure map of $M$.  Then
\[
RM=\{x\in M |\mu ^{*}(x)=\alpha (y) \text{ for some } y\}.  
\]
We refer to $RM$ as the largest subcomodule of $M$.   One can then
define the product of a set of comodules $\{M_{i} \}$ to be the largest
subcomodule of the $A$-module product.  
\end{remark}

Building on the remark above, note that, for a set of comodules $\{M_{i}
\}$, there is the natural commutative diagram of $A$-modules below.
\[
\begin{CD}
\prod_{i}^{\Gamma} M_{i} @>>> \Gamma \otimes \prod_{i} M_{i} @>>> \Gamma
\otimes \prod_{i} N_{i} \\
@. @V\alpha VV @V\alpha VV \\
\prod_{i} M_{i} @>>\prod_{i}\psi > \prod_{i} (\Gamma \otimes M_{i})
@>>\prod_{i}\psi p_{i} > \prod_{i} (\Gamma \otimes N_{i}).  
\end{CD}
\]
This means that there is a natural induced map of $A$-modules
$\prod_{i}^{\Gamma} M_{i} \xrightarrow{} \prod_{i}M_{i}$.  This map is
injective when $\alpha$ is injective, which is certainly true if
$\Gamma$ is projective over $A$.  It is an isomorphism when $\alpha$ is
so, which is true if $\Gamma$ is finitely generated and projective over
$A$.

Since the product is right adjoint to the exact diagonal functor, the
product is left exact.  But it need not be exact in general.  Indeed,
let $A=\mathbf{Q}$ and $\Gamma =A[x]$, thought of as a primitively
generated Hopf algebra over $A$.  Let $X_{n}=\Gamma /(x^{n})$ for $n\geq
1$, and let $Y_{n}=A$.  There is a surjection $X_{n}\xrightarrow{}Y_{n}$
that sends $x^{n-1}$ to $1$ and every other power of $x$ to $0$.  But
one can check that $\prod^{\Gamma}Y_{n}\cong \prod_{n}Y_{n}$, and that
there is no element of $\prod^{\Gamma}X_{n}$ that hits
$(1,1,\dotsc,1,\dotsc)$.  Indeed, $\prod^{\Gamma} X_{n}$ consists of
those elements $(f_{1},f_{2},\dotsc)$ of $\prod X_{n}$ such that the
degrees of $f_{i}$ are bounded.
 
We can also use this technique of constructing right adjoints to prove
the following proposition.  

\begin{proposition}\label{prop-comodules-adjunction}
Suppose $\Phi\mathcolon \hopfalg \xrightarrow{} \otherhopfalg$ is a map
of Hopf algebroids.  Then the functor $\Phi_{*}\mathcolon \Gamma \comod
\xrightarrow{}\Sigma \comod$ has a right adjoint $\Phi^{*}$.
\end{proposition}

\begin{proof}
An adjointness argument shows that we must define 
\[
\Phi^{*}(\Sigma \otimes_{B}N)=\Gamma \otimes N
\]
when $N$ is a $B$-module.  Given a comodule map $f\mathcolon \Sigma
\otimes_{B}N\xrightarrow{}\Sigma \otimes_{B}N'$, we define
$\Phi^{*}(f)\mathcolon \Gamma \otimes N\xrightarrow{}\Gamma \otimes N'$
as the following composite.
\[
\Gamma \otimes N\xrightarrow{\Delta \otimes 1} \Gamma \otimes \Gamma
\otimes N \xrightarrow{1\otimes \alpha} \Gamma \otimes \Sigma
\otimes_{B}N \xrightarrow{1\otimes f} \Gamma \otimes \Sigma \otimes_{B}
N' \xrightarrow{1\otimes \epsilon \otimes 1} \Gamma \otimes N'.  
\]
Here the map $\alpha$ is defined by $\alpha (x\otimes
n)=\Phi_{1}(x)\otimes n$.  We leave it to the reader to check that this
definition is functorial, so that we have defined $\Phi^{*}$ on the full
subcategory of extended comodules.  

As usual, given an arbitrary $\Sigma$-comodule $N$, we write $N$ as the
kernel of 
\[
\psi p\mathcolon \Sigma \otimes_{B} N\xrightarrow{}\Sigma
\otimes_{B} N',
\]
where $p\mathcolon \Sigma \otimes_{B} N\xrightarrow{}N'$ is the cokernel
of $\psi$.  We then define $\Phi^{*}(N)$ as the kernel of
\[
\Gamma \otimes N \xrightarrow{\Phi^{*}(\psi p)} \Gamma \otimes N'.
\]
We leave to the reader the check that $\Phi^{*}$ is right adjoint to
$\Phi_{*}$.  
\end{proof}

\subsection{Duality and finite presentation}\label{subsec-comod-duality}

We now show that $\Gamma $-comod is in fact a closed symmetric monoidal
category, and we characterize the dualizable comodules.  

\begin{theorem}\label{thm-comod-closed}
If $\hopfalg$ is a flat Hopf algebroid, then the category $\Gamma\comod$
is closed symmetric monoidal.  Furthermore, the closed structure
$F(M,N)$ is left exact in $N$ and right exact in $M$.
\end{theorem}
 
\begin{proof}
An adjointness argument shows that we must define 
\[
F(M,\Gamma \otimes N)=\Gamma \otimes \Hom_{A}(M,N).
\]
Suppose we have a map $f\mathcolon \Gamma \otimes N\xrightarrow{f}\Gamma
\otimes N'$ of extended comodules.  We need to define the map
$F(M,f)\mathcolon \Gamma \otimes \Hom_{A}(M,N)\xrightarrow{}\Gamma
\otimes \Hom_{A}(M,N')$.  This map will be adjoint to a map
\[
\Gamma \otimes \Hom_{A}(M,N) \xrightarrow{} \Hom_{A}(M,N')
\]
of $A$-modules, which will in turn be adjoint to the composite 
\[
\Gamma \otimes \Hom_{A}(M,N)\otimes M\xrightarrow{1\otimes \ev} \Gamma
\otimes N\xrightarrow{f} \Gamma \otimes N'\xrightarrow{\epsilon \otimes
1} N'.  
\]
We leave to the reader the check that this definition is functorial, so
that we have defined $F(M,-)$ on the full subcategory of extended
comodules.  We also leave to the reader the check that there is a
natural isomorphism 
\[
\Gamma \comod (L,F(M,\Gamma \otimes N)) \cong \Gamma \comod (L\smash 
M,\Gamma \otimes N),
\]
where naturality refers to an arbitrary map of extended comodules
$\Gamma \otimes N\xrightarrow{}\Gamma \otimes N'$.  

We then have no choice but to define $F(M,N)$ as the kernel of
\[
\Gamma \otimes \Hom (M,N) \xrightarrow{F(M,\psi p)}\Gamma \otimes \Hom (M,N')
\]
where $p\mathcolon \Gamma \otimes N\xrightarrow{}N'$ is the cokernel of
$\psi$.  The necessary adjunction isomorphism follows immediately.   

Since $F(M,-)$ is right adjoint to the right exact functor $M\otimes -$,
it is left exact.  Now suppose we have a right exact sequence
\[
M'\xrightarrow{}M\xrightarrow{}M''\xrightarrow{} 0.
\]
Then we have a right exact sequence 
\[
L\smash  M' \xrightarrow{} L\smash M \xrightarrow{} L\smash M''
\xrightarrow{}0, 
\]
and so a left exact sequence of abelian groups 
\[
0\xrightarrow{} \Gamma \comod (L\smash M'',N) \xrightarrow{} \Gamma
\comod (L\smash M,N) \xrightarrow{} \Gamma \comod (L\smash M',N).
\]
Applying adjointness, we find that $F(M'',N)$ has the universal property
characterizing the kernel of $F(M,N)\xrightarrow{}F(M',N)$.  
\end{proof}

It would be nice to have a better understanding of $F(M,N)$.  The
following proposition is helpful.  

\begin{proposition}\label{prop-comodules-function}
Suppose $\hopfalg$ is a flat Hopf algebroid, and $M$ and $N$ are
$\Gamma$-comodules.  
\begin{enumerate}
\item [(a)] There is a natural map
$F(M,N)\xrightarrow{\tau_{MN}}\Hom_{A}(M,N)$ of $A$-modules.
\item [(b)] If $M$ is finitely generated over $A$, then $\tau_{MN}$ is
injective.
\item [(c)] If $M$ is finitely presented over $A$, then $\tau_{MN}$ is
an isomorphism.  
\end{enumerate}
\end{proposition}

\begin{proof}
Consider the natural diagram below.  
\[
\begin{CD}
F(M,N) @>>> \Gamma \otimes \Hom_{A} (M,N) @>>> \Gamma \otimes \Hom_{A}
(M,N') \\ 
@. @VVV @VVV \\
\Hom_{A} (M,N) @>>> \Hom_{A} (M,\Gamma \otimes N) @>>> \Hom_{A}(M,\Gamma
\otimes N')
\end{CD}
\]
The vertical arrows take $x\otimes f$ to the map that takes $m$ to
$x\otimes f(m)$.  It is not obvious that this diagram is commutative,
but a careful diagram chase shows that it is.  The rows both express
their left-hand entry as a kernel, the first row by definition, and the
second row by applying $\Hom_{A}(M,-)$ to diagram~\ref{eq-comod}.  Thus,
there is a natural induced map $F(M,N)\xrightarrow{}\Hom_{A}(M,N)$,
proving part~(a).

Parts~(b) and~(c) will follow if we can show that the vertical maps are
injections when $M$ is finitely generated and isomorphisms when $M$ is
finitely presented.  Since $\Gamma$ is flat as a right $A$-module, we
can write $\Gamma =\colim C_{i}$, where the $C_{i}$ are finitely
generated projective $A$-modules.  The natural map 
\[
C_{i}\otimes \Hom_{A}(M,N)\xrightarrow{}\Hom_{A}(M,C_{i}\otimes N)
\]
is therefore an isomorphism.  Hence, the map
\[
\Gamma \otimes \Hom_{A} (M,N)\cong \colim_{i} C_{i}\otimes \Hom_{A} (M,N)
\xrightarrow{} \colim_{i} \Hom_{A}(M,C_{i}\otimes N)
\]
is also an isomorphism.  When $M$ is finitely generated over $A$, the
map 
\[
\colim_{i} \Hom_{A}(M, C_{i}\otimes N)\xrightarrow{}\Hom_{A}(M, \Gamma
\otimes N)
\]
is injective; when $M$ is finitely presented over $A$, it is an
isomorphism.  Parts~(b) and~(c) follow.
\end{proof}

We now recall that an object $M$ in a cocomplete category $\cat{C}$ is
called \emph{$\lambda$-presented}, for a regular cardinal $\lambda$, if
$\cat{C}(M,-)$ commutes with $\lambda$-filtered colimits
(See~\cite[Section~6.4]{borceux-1}).  When $\lambda =\omega $, we get
the usual notion of a \emph{finitely presented} object.  An $A$-module
$M$ is $\lambda $-presented if and only if it is a quotient of a map of
free modules, each of which has rank $<\lambda $.  

\begin{proposition}\label{prop-comodules-presented}
Suppose $\hopfalg$ is a flat Hopf algebroid, $M$ is a $\Gamma
$-comodule, and $\lambda $ is a regular cardinal.  Consider the
following three statements.  
\begin{enumerate}
\item [(a)] $M$ is $\lambda $-presented as a $\Gamma $-comodule. 
\item [(b)] $M$ is $\lambda $-presented as an $A$-module. 
\item [(c)] The functor $F(M,-)$ commutes with $\lambda $-filtered
colimits.  
\end{enumerate}
Then $(a)$ and $(b)$ are equivalent, and $(c)$ implies $(a)$.  In
particular, $M$ is $\kappa $-presented for some $\kappa $.
\end{proposition}

\begin{proof}
We first show that (a) implies (b).  So suppose that $M$ is $\lambda
$-presented as a $\Gamma $-comodule, and we have a $\lambda $-filtered
diagram of $A$-modules $N_{i}$.  Then 
\begin{gather*}
\colim \Hom _{A}(M, N_{i}) \cong \colim \Gamma \comod (M, \Gamma \otimes
N_{i}) \\
\cong \Gamma \comod (M, \colim (\Gamma \otimes N_{i})) \cong \Gamma
\comod (M, \Gamma \otimes \colim N_{i}) \\
\cong \Hom _{A}(M, \colim N_{i}).  
\end{gather*}

We now show that (b) implies (a).  Suppose that $M$ is $\lambda
$-presented as an $A$-module, and we have a $\lambda $-filtered diagram
of comodules $N_{i}$.  We must show that the map
\[
\colim \Gamma \comod (M,N_{i}) \xrightarrow{} \Gamma \comod (M, \colim
N_{i})
\]
is an isomorphism.  It is obviously injective, since the forgetful
functor to $A$-modules is faithful, and $M$ is $\lambda$-presented as an
$A$-module.  On the other hand, suppose we have a map $f\mathcolon
M\xrightarrow{}\colim N_{i}$ of comodules.  As a map of $A$-modules $f$
factors through some map $g\mathcolon M\xrightarrow{}N_{i}$ for
some $i$.  The difficulty is that $g$ may not be a map of
comodules, since $\psi g$ may not be equal to $(1\otimes g)\psi$.  But
they are equal as maps to 
\[
\Gamma \otimes \colim N_{j}\cong \colim (\Gamma \otimes N_{j}),
\]
so they must be equal in some $\Gamma \otimes N_{j}$.
It follows that the composite
\[
M \xrightarrow{g} N_{i} \xrightarrow{} N_{j}
\]
is the desired factorization of $f$.  

We now show that (c) implies (a).  So suppose that $F(M,-)$
commutes with $\lambda $-filtered colimits, and we have a $\lambda
$-filtered diagram of comodules $N_{i}$.  Then 
\begin{gather*}
\colim \Gamma \comod (M,N_{i}) \cong \colim \Gamma \comod (A,
F(M,N_{i})) \\
\cong \Gamma \comod (A, \colim F(M,N_{i})) \cong \Gamma \comod (A,
F(M, \colim N_{i})) \\
\cong \Gamma \comod (M, \colim N_{i}),
\end{gather*}
where the second isomorphism holds because $A$ is finitely presented.  
\end{proof}

In any closed symmetric monoidal category with unit $A$, we define
$DM=F(M,A)$.  There is always a natural map $DM\smash
N\xrightarrow{}F(M,N)$.  When this map is an isomorphism for all $N$,
$M$ is called \emph{strongly dualizable}, which we generally abbreviate
to \emph{dualizable}.  The author does not know to whom this concept is
due; perhaps Puppe~\cite{puppe-duality}.  An excellent reference
is~\cite[Chapter~III]{lewis-may-steinberger}, and the basic properties
of dualizable objects are summarized
in~\cite[Theorem~A.2.5]{hovey-axiomatic}.

\begin{proposition}\label{prop-comodules-dualizable}
Suppose $\hopfalg$ is a flat Hopf algebroid.  Then a $\Gamma $-comodule
$M$ is dualizable in $\Gamma \comod $ if and only if $M$ is finitely
generated and projective as an $A$-module.
\end{proposition}

\begin{proof}
First suppose that $M$ is finitely generated and projective as an
$A$-module.  We need to check that the map $F(M,A)\smash
N\xrightarrow{}F(M,N)$ is an isomorphism.  But $F(M,A)\cong
\Hom_{A}(M,A)$ and $F(M,N)\cong \Hom_{A}(M,N)$ by
Proposition~\ref{prop-comodules-function}.  It is well known and easy to
check that the natural map $\Hom_{A}(M,A)\otimes
N\xrightarrow{}\Hom_{A}(M,N)$ is an isomorphism when $M$ is a finitely
generated projective.

Now suppose that $M$ is dualizable.  Then the functor $F(M,-)$ commutes
with colimits, since it is isomorphic to $F(M,A)\smash (-)$.
Proposition~\ref{prop-comodules-presented} then implies that $M$ is
finitely presented as an $A$-module.  We now show that $M$ must in fact
be projective over $A$.  Indeed, the functor $F(M,-)$ is always left
exact, and because $M$ is dualizable, $F(M,-)\cong F(M,A)\smash (-)$ is
also right exact.  Hence $F(M,-)$ is an exact functor on the category of
$\Gamma $-comodules.  But Proposition~\ref{prop-comodules-function}
tells us that $F(M,-)\cong \Hom _{A}(M,-)$ since $M$ is finitely
presented over $A$.  Now suppose $E$ is an exact sequence of
$A$-modules.  Then $\Hom_{A}(M, \Gamma \otimes E)$ is again exact.  But,
since $M$ is finitely presented, 
\[
\Hom_{A}(M, \Gamma \otimes E)\cong \Gamma \otimes \Hom_{A}(M,E)
\]
by the argument of Proposition~\ref{prop-comodules-function}.  Since
$\Gamma$ is faithfully flat over $A$, we conclude that $\Hom_{A}(M,E)$
is exact, so $M$ is projective over $A$.
\end{proof}

\subsection{Generators and Adams Hopf algebroids}\label{subsec-Adams}

We have just seen that the category of $\Gamma $-comodules has many good
properties when $(A, \Gamma )$ is a flat Hopf algebroid.  But those
properties are still not enough for us, because we need a good set of
generators for the category of $\Gamma $-comodules.  Recall that a set
of objects $\cat{G}$ in an abelian category $\cat{C}$ is said to
\emph{generate} $\cat{C}$ when $\cat{C}(P,f)=0$ for all $P\in \cat{G}$
implies that $f=0$.

For much of the sequel, we will require that the dualizable comodules
generate the category of $\Gamma $-comodules.  The main advantage of
this hypothesis is the following proposition.  

\begin{proposition}\label{prop-generators}
Suppose $\hopfalg $ is a flat Hopf algebroid for which the dualizable
comodules generate the category of $\Gamma $-comodules.  Then the
category of $\Gamma $-comodules is a locally finitely presentable
Grothendieck category.  Furthermore\uc
\begin{enumerate}
\item [(a)] Every comodule is a quotient of a direct sum of dualizable
comodules.  
\item [(b)] If $M$ is a comodule and $x\in M$, then there is a
dualizable comodule $P$ and a map $P\xrightarrow{}M$ of comodules whose
image contains $x$.  
\item [(c)] Every comodule that is finitely generated over $A$ is a
quotient of a dualizable comodule.
\item [(d)] Every comodule is the union of its subcomodules that are
finitely generated over $A$.  
\item [(e)] Every comodule is a filtered colimit of finitely presented
comodules. 
\end{enumerate}
\end{proposition}

Most of these facts are true in a general locally finitely
presentable Grothendieck category; see~\cite{stenstrom} for details.
We therefore give only a sketch of the proof.  

\begin{proof}
The first statement is true because there is only a set of isomorphism
classes of dualizable comodules, and dualizable comodules are finitely
presented.  For part~(a), let $\cat{G}$ denote a set containing one
element from each isomorphism class of dualizable comodules, and let $T$
be the set of all maps $f$ with $\dom f\in \cat{G}$ and $\codom f=M$.
Consider the map
\[
\alpha \mathcolon \bigoplus_{f\in T} \dom f\xrightarrow{}M,
\]
and let $\beta $ denote the cokernel of this map.  Then $\Gamma \comod
(P,\beta )=0$ for all dualizable $P$, so $\beta =0$.  Hence $\alpha $ is
surjective.  

Part~(b) is an immediate corollary of part~(a), and part~(c) and
part~(d) follow easily from part~(b).  For part~(e), we choose a small
skeleton $\cat{F}$ of the category of finitely presented comodules and
consider the category $\cat{F}/M$ consisting of all maps from an element
of $\cat{F}$ to $M$.  There is an obvious map 
\[
\colim _{f\in \cat{F}/M} \dom f \xrightarrow{} M.
\]
One can readily verify that $\cat{F}$ is filtered, and that this map is
a monomorphism, for any flat Hopf algebroid $(A, \Gamma )$.  If the
dualizable comodules generate, then it is an epimorphism by part~(b).
\end{proof}

We also have the following corollary, which answers the question left
open by Proposition~\ref{prop-comodules-presented}.  

\begin{corollary}\label{cor-comodules-presented}
Suppose $\hopfalg $ is a flat Hopf algebroid for which the dualizable
comodules generate the category of $\Gamma $-comodules, $M$ is a $\Gamma
$-comodule, and $\lambda $ is a regular cardinal.  If $M$ is $\lambda
$-presented, then $F(M, -)$ commutes with $\lambda $-filtered colimits.
\end{corollary}

\begin{proof}
Suppose we have a $\lambda $-filtered system of comodules $N_{i}$. We
need to show that 
\[
\colim F(M,N_{i})\xrightarrow{\alpha }F(M, \colim N_{i})
\]
is an isomorphism.  Because the dualizable comodules generate, it
suffices to show that $\Gamma \comod (P,f)$ is an isomorphism for all
dualizable comodules $P$.  Since dualizable comodules are in particular
finitely presented, we have
\begin{gather*}
\Gamma \comod (P, \colim F(M,N_{i})) \cong \colim \Gamma \comod (P,
F(M,N_{i})) \\
\cong \colim \Gamma \comod (P\smash M, N_{i}) \cong \colim \Gamma
\comod (M, DP \smash N_{i}) \\
\cong \Gamma \comod (M, \colim (DP \smash N_{i})) \cong \Gamma \comod
(M, DP \smash \colim N_{i}) \\
\cong \Gamma \comod (M\smash P, \colim N_{i}) \cong \Gamma \comod (P,
F(M, \colim N_{i})).  
\end{gather*}
\end{proof}

We now owe the reader some examples of flat Hopf algebroids for which
the dualizable comodules generate the category of $\Gamma $-comodules.
We learned the following definition
from~\cite{goerss-hopkins-comodules}, but it is implicit
in~\cite[Section~III.13]{adams-blue}.

\begin{definition}\label{defn-Adams}
A Hopf algebroid $\hopfalg $ is said to be an \emph{Adams Hopf
algebroid} when $\Gamma $ is the colimit of a filtered system of
comodules $\Gamma _{i}$, where $\Gamma _{i}$ is finitely generated and
projective over $A$.
\end{definition}

In particular, any Adams Hopf algebroid is flat, since the colimit of
projective modules is flat.  Thus the $\Gamma _{i}$ are dualizable
comodules.  The following proposition is a restatement
of~\cite[Lemma~3.4]{goerss-hopkins-comodules}.

\begin{proposition}\label{prop-amenable-Adams}
If $\hopfalg $ is an Adams Hopf algebroid, then the dualizable comodules
generate the category of $\Gamma $-comodules.  
\end{proposition}

\begin{proof}
Suppose $\hopfalg $ is Adams, and $M$ is a $\Gamma $-comodule. Then we
have:
\begin{gather*}
M \cong \Hom_{A}(A,M) \cong \Gamma \comod (A,\Gamma \otimes M) \\
\cong \Gamma \comod (A, \Gamma \smash M) \cong 
\Gamma \comod (A,\colim \Gamma_{i}\smash M) \\
\cong \colim \Gamma \comod (A,\Gamma_{i}\smash M) \cong \colim \Gamma
\comod (D\Gamma_{i},M).  
\end{gather*}
The result follows.  
\end{proof}
  
We now give some examples of Adams Hopf algebroids.  Most of the ones we
are interested in come from algebraic topology.  Recall the notion of
minimal weak colimit from~\cite[Section~2.2]{hovey-axiomatic}.

\begin{definition}\label{defn-top-flat}
A ring spectrum $R$ is called \emph{topologically flat} if $R$ is the
minimal weak colimit of a filtered diagram of finite spectra $X_{i}$
such that $R_{*}X_{i}$ is a finitely generated projective
$R_{*}$-module.
\end{definition}

This definition is based on~\cite[Condition~III.13.3]{adams-blue}.  

\begin{lemma}\label{lem-top-flat}
Suppose $R$ is a ring spectrum that is topologically flat and such that
$R_{*}R$ is commutative.  Then $(R_{*},R_{*}R)$ is an Adams Hopf
algebroid.  
\end{lemma}

\begin{proof}
Write $R$ as the minimal weak colimit of the $X_{i}$.  Then
$R_{*}R=\colim R_{*}X_{i}$.  In particular, $R_{*}R$ is flat over
$R_{*}$ (and this is the reason for the term ``topologically flat'') and
satisfies the Adams condition.  
\end{proof}

The reason for the hypothesis that $R_{*}R$ be commutative is that there
could well be non-commutative ring spectra $R$, such as Morava
$K$-theory $K(n)$ at the prime $2$, where $R_{*}R$ is nevertheless
commutative.

Adams gave several examples of topologically flat ring spectra
in~\cite[Proposition~III.13.4]{adams-blue}, which we restate here.

\begin{theorem}\label{thm-Adams}
The ring spectra $MU$, $MSp$, $K$, $KO$, $H\Fp $, and $K(n)$ are
topologically flat.  
\end{theorem}

Adams did not of course consider $K(n)$, since it had not been
discovered yet, but his proof for $H\Fp $ works for any field spectrum.  

We can add another case to this list as well.  Recall that $BP$ is the
Brown-Peterson spectrum, and $BP_{*}\cong \Zp[v_{1}, v_{2},\dots ]$.
Given an invariant regular sequence $J=(p^{i_{0}}, v_{1}^{i_{1}},\dots
,v_{k-1}^{i_{k-1}})$ in $BP_{*}$, there is a spectrum $BPJ$ with
$BPJ_{*}\cong BP_{*}/J$ studied in~\cite{johnson-yosimura}.

\begin{proposition}\label{prop-BPJ-Adams}
Let $J=(p^{i_{0}}, v_{1}^{i_{1}},\dots ,v_{k-1}^{i_{k-1}})$ be an
invariant regular sequence of length $k$ in $BP_{*}$.  Then $BPJ$ is
topologically flat.  
\end{proposition}

\begin{proof}
Write $BP$ as a minimal weak colimit of spectra $X_{\alpha }$, where
$X_{\alpha }$ is a finite spectrum with cells in only even degrees. 
Then $BPJ\smash BP$ is the minimal weak colimit of the $BPJ\smash
X_{\alpha }$.  On the other hand, we claim that $BPJ\smash BPJ$ is a
wedge of $2^{k}$ copies of $BPJ\smash BP$.  Indeed, $BPJ_{*}BPJ$ is a
free module over $BPJ_{*}BP$, and so we choose generators for the free
module and use them to construct the desired splitting.  

Hence $BPJ\smash BPJ$ is the minimal weak colimit of $BPJ \smash
Y_{\alpha }$, where $Y_{\alpha }$ is a finite wedge of copies of
$X_{\alpha }$.  Since each $BPJ_{*}(Y_{\alpha })$ is a free module, this
completes the proof. 
\end{proof}

The following theorem, generalizing Proposition~2.12
of~\cite{hovey-strickland}, gives us many other examples of
topologically flat ring spectra.   Recall that, if $R$ is a ring
spectrum and $E$ is an $R$-module spectrum, then $E$ is said to be
\emph{Landweber exact} over $R$ if the natural map 
\[
E_{*}\otimes _{R_{*}} R_{*}X \xrightarrow{} R_{*}X
\]
is an isomorphism for all spectra $X$. 
 
\begin{theorem}\label{thm-Adams-Landweber}
Suppose $R$ is a topologically flat ring spectrum, and $E$ is a
Landweber exact $R$-module spectrum.  Then $E$ is topologically flat.  
\end{theorem}

\begin{proof}
The proof is much like that of Proposition~2.12
of~\cite{hovey-strickland}.  Write $R$ as the minimal weak colimit of
finite spectra $X_{i}$ such that $R_{*}X_{i}$ is finitely generated and
projective over $R_{*}$.  Then $E_{*}X_{i}$ is finitely generated and
projective over $E_{*}$.  We show that $E$ is the filtered colimit of a
diagram of finite wedges of suspensions of $X_{i}$.  To do so, consider
the category $\cat{F}/E$ of all maps of finite spectra to $E$, and
consider the full subcategory $\cat{F}'/E$ of such maps whose domain is
a (variable) finite wedge of suspensions of the $X_{j}$.  We claim that
this is cofinal in $\cat{F}/E$.  Since $E$ is the minimal weak colimit
of the obvious functor from $\cat{F}/E$ to spectra, it will follow that
$E$ is the minimal weak colimit of the restriction of this functor to
$\cat{F}'/E$.  

To show that $\cat{F}'/E$ is cofinal in $\cat{F}/E$, it suffices to show
that any map $f$ from a finite spectrum $Z$ to $E$ factors through a
finite wedge of suspensions of the $X_{i}$.  By Spanier-Whitehead
duality
\[
E^{*}(Z) \cong E^{*} \otimes _{R^{*}} R^{*}Z.
\]
We can thus write $f=\sum _{i=1}^{m} b_{i}\otimes c_{i}$.  Because $R$
is the minimal weak colimit of the $X_{j}$, each map
$c_{i}$ has a factorization
\[
c_{i} = (Z \xrightarrow{g_{i}}\Sigma ^{-|c_{i}|}X_{i}
\xrightarrow{e_{i}}\Sigma ^{-|c_{i}|}R).
\]
Let $Y=\bigvee_{i=1}^{m} \Sigma ^{-|c_{i}|}X_{i}$, let $g\mathcolon
Z\xrightarrow{}Y$ be the map with components $g_{i}$, and let
$h\mathcolon Y\xrightarrow{}E$ be the map with components $b_{i}\otimes
e_{i}\in E^{*}\otimes _{R^{*}}R^{*}X_{i}\cong E^{*}(X_{i})$.  This gives
the desired factorization.  
\end{proof}

Of course, there are algebraic examples of Adams Hopf algebroids as
well.  

\begin{proposition}\label{prop-algebra-Adams}
Suppose $\Gamma $ is a Hopf algebra over a field $k$.  Then $(k,\Gamma )$ is
an Adams Hopf algebroid.  
\end{proposition}

\begin{proof}
By Lemma~9.5.3 of~\cite{hovey-axiomatic}, every $\Gamma$-comodule is the
filtered colimit of its finite-dimensional sub-comodules.  In
particular, this is true for $\Gamma$ itself.  
\end{proof}

\begin{proposition}\label{prop-Adams-quotients}
Suppose $\hopfalg $ is an Adams Hopf algebroid, $I$ is an invariant ideal in
$A$, and $v$ is a primitive element in $A$.  Then $(A/I, \Gamma /I)$ and
$(v^{-1}A, v^{-1}\Gamma )$ are Adams Hopf algebroids. 
\end{proposition}

\begin{proof}
Suppose that $\Gamma \cong \colim \Gamma _{j}$, where
each $\Gamma _{j}$ is finitely generated and projective over $A$.  Then
$\Gamma /I=\colim  \Gamma _{j}/I$ and $\Gamma _{j}/I$ is finitely
generated and projective over $A/I$.  Similarly, $v^{-1}\Gamma =\colim
v^{-1}\Gamma _{j}$.  
\end{proof}

Despite all these examples of Adams Hopf algebroids, there is a
theoretical difficuly with the notion.  

\begin{question}\label{quest-Adams}
Suppose $\Phi \mathcolon \hopfalg \xrightarrow{}\otherhopfalg $ is a
weak equivalence of Hopf algebroids.  Is it true that $\hopfalg $ is
Adams if and only if $\otherhopfalg $ is Adams?  
\end{question}

Note that if $\Phi $ is a weak equivalence, then the dualizable $\Gamma
$-comodules generate if and only if the dualizable $\Sigma $-comodules
generate.

\section{The projective model structure}\label{sec-model-proj}

In this section, we establish a preliminary model structure on
$\Ch{\Gamma}$, the category of unbounded chain complexes of
$\Gamma$-comodules.

\subsection{Construction and basic properties}\label{subsec-proj-basic}

We recall the results of~\cite{hovey-christensen-relative}.  Beginning
with a set of objects $\cat{S}$ in a cocomplete abelian category
$\cat{A}$, there is a projective class $(\cat{P},\cat{E})$, where
$\cat{E}$ consists of all maps $f$ such that $\cat{A}(P,f)$ is onto for
all $P$ in $\cat{S}$, and $\cat{P}$ consists of all retracts of direct
sums of elements of $\cat{S}$.  The elements of $\cat{P}$ are called
\emph{relative projectives}, and the maps of $\cat{E}$ are called
\emph{relative epimorphisms}.  This is
~\cite[Lemma~1.5]{hovey-christensen-relative}, but it is also easy to
see.

The main result of~\cite{hovey-christensen-relative} associates a model
structure on $\Ch{\cat{A}}$, the category of unbounded chain complexes
in $\cat{A}$, to a projective class $(\cat{P},\cat{E})$, given some
hypotheses.  We recall that a chain map $\phi$ is a fibration in this
model structure when $\cat{A}(P,\phi)$ is a degreewise surjection for
all $P\in \cat{P}$, and a weak equivalence when $\cat{A}(P,\phi)$ is a
homology isomorphism for all $P\in \cat{P}$.

The hypothesis needed is that functorial cofibrant replacements exist.
This is automatic,
by~\cite[Proposition~4.2]{hovey-christensen-relative}, when $\cat{A}$ is
complete and cocomplete, there are enough $\kappa$-small
$\cat{P}$-projectives for some cardinal $\kappa$, and functorial
$\cat{P}$-resolutions exist.  When $\cat{P}$ is generated by a set
$\cat{S}$ as above, then functorial $\cat{P}$-resolutions obviously
exist, since there is a functorial $\cat{P}$-epic
\[
\bigoplus_{P\in \cat{S}} \bigoplus_{f\in \cat{A}(P,M)} P \xrightarrow{} M
\]
for any $M\in \cat{A}$.  When each object of $\cat{S}$ is $\lambda
$-small for some $\lambda $, then there are enough $\kappa$-small
$\cat{P}$-projectives (take $\kappa$ to be the supremum of the
$\lambda$'s).

Now, if $\cat{A}$ happens to be a Grothendieck abelian category, then it
is automatically complete and cocomplete, and every object in $\cat{A}$
is $\kappa$-presented, and so \emph{a fortiori} $\kappa$-small, for some
$\kappa$.  This latter statement is an immediate corollary of the fact
that Grothendieck abelian categories are locally
presentable~\cite[Proposition~3.10]{beke}, but a direct proof can be
found in the Appendix to~\cite{hovey-sheaves}.

We thus have the following result, which was inexplicably not stated
in~\cite{hovey-christensen-relative}.

\begin{theorem}\label{thm-christensen-relative}
Suppose $\cat{A}$ is a Grothendieck abelian category, and $\cat{S}$ is a
set of objects in $\cat{A}$.  Then there is a model structure on
$\Ch{\cat{A}}$ in which the fibrations are the maps $\phi$ such that
$\cat{A}(P,\phi )$ is a surjection for all $P\in \cat{S}$ and the weak
equivalences are the maps $\phi$ such that $\cat{A}(P,\phi)$ is a
homology isomorphism for all $P\in \cat{S}$.  More generally, this model
structure exists when $\cat{A}$ is complete and cocomplete, but not
necessarily Grothendieck, as long as every object of $\cat{S}$ is
$\kappa$-small for some $\kappa$.
\end{theorem}

Now we return to the case at hand, when $\cat{A}$ is the category of
$\Gamma$-comodules and $\hopfalg$ is a Hopf algebroid.  In the light of
the results of Section~\ref{subsec-Adams}, we should take $\cat{S}$ to
be the set of dualizable $\Gamma $-comodules.  

\begin{definition}\label{defn-projective}
Suppose $\hopfalg$ is a flat Hopf algebroid.  Let $\cat{S}$ be a set
containing one comodule from each isomorphism class of dualizable
$\Gamma$-comodules.  We refer to the retracts of direct sums of elements
of $\cat{S}$ as \emph{relatively projective} comodules, and to the maps
$f$ of comodules such that $\Gamma \comod (P,f)$ is surjective for all
$P\in \cat{S}$ as \emph{relative epimorphisms}.  The resulting model
structure on $\Ch{\Gamma}$ obtained from
Theorem~\ref{thm-christensen-relative} is called the \emph{projective
model structure}.  Thus, a map $\phi$ is a \emph{projective fibration}
if $\phi $ is a degreewise relative epimorphism, and $\phi$ is a
\emph{projective equivalence} if $\Gamma \comod (P,\phi)$ is a homology
isomorphism for all $P\in \cat{S}$.  The map $\phi$ is a
\emph{projective cofibration}, or simply a \emph{cofibration}, if $\phi$
has the \llp all projective trivial fibrations.  We refer to an chain
complex $F$ as \emph{projectively trivial} if $0\xrightarrow{}F$ is a
projective equivalence.
\end{definition}

Goerss and Hopkins~\cite{goerss-hopkins-comodules} put a model structure
on the category of nonnegatively graded chain complexes over an Adams
Hopf algebroid.  Their model structure gave us the idea for the
projective model structure, but it is not the same, as they took
$\cat{S}$ to be the set of all the $D\Gamma _{i}$ (under the assumption
that $\Gamma =\colim \Gamma _{i}$).  One obvious advantage of our
definition is that the dualizable comodules are canonically attached to
the symmetric monoidal category of $\Gamma $-comodules, while the
$D\Gamma _{i}$ are not.

We point out that we do not need $\hopfalg$ to be an Adams Hopf
algebroid, or even for the dualizable comodules to generate, for the
projective model structure to exist.  Also note that every relatively
projective comodule is projective as an $A$-module, but we do not know
if the converse holds.

Note also that because the elements of $\cat{S}$ are finitely presented
in the category $\Gamma \comod$ (see
Proposition~\ref{prop-comodules-presented}), filtered colimits of
projective equivalences (resp. projective fibrations) are again
projective equivalences (resp. projective fibrations).

We then have the following theorem describing some of the properties of
the projective model structure.  
 
\begin{theorem}\label{thm-projective-model}
Suppose $\hopfalg$ is a flat Hopf algebroid.  Then the projective model
structure on $\Ch{\Gamma}$ is proper, finitely generated, stable, and
symmetric monoidal.  A map $\phi$ is a cofibration if and only if it is
a degreewise split monomorphism whose cokernel is cofibrant. A chain
complex $X$ is cofibrant if and only if it is a retract of a colimit of
complexes
\[
X_{0} \xrightarrow{} X_{1} \xrightarrow{}\dotsb X_{\alpha } \dotsb 
\]
where each $X_{\alpha}\xrightarrow{}X_{\alpha +1}$ is a degreewise split
monomorphism whose cokernel is a complex of relative projectives with no
differential.  The homotopy relation between cofibrant objects is the
usual chain homotopy relation.
\end{theorem}

\begin{proof}
This all follows from the results of~\cite{hovey-christensen-relative}.
The characterization of cofibrant objects follows from Corollary~4.4
of~\cite{hovey-christensen-relative}, and the characterization of
cofibrations follows from Proposition~2.5
of~\cite{hovey-christensen-relative}.  The fact that the model structure
is proper is~\cite[Proposition~2.18]{hovey-christensen-relative}, and
stability is~\cite[Corollary~2.17]{hovey-christensen-relative} and
obvious.  The fact that homotopy is the usual chain homotopy
is~\cite[Lemma~2.13]{hovey-christensen-relative}.

The generating cofibrations for the projective model structure are
$S^{n-1}P\xrightarrow{}D^{n}P$ for $P\in \cat{S}$, and the generating
trivial cofibrations are $0\xrightarrow{}D^{n}P$ for $P\in \cat{S}$.
Here $S^{n-1}P$ denotes the complex which is $P$ in degree $n-1$ and
zero elsewhere, and $D^{n}P$ denotes the complex which is $P$ in degrees
$n$ and $n-1$ and $0$ elsewhere.  This is proved in Section~5
of~\cite{hovey-christensen-relative}.  Each of $S^{n-1}P$, $0$, and
$D^{n}P$ are finitely presented, so the projective model structure is
finitely generated.

Finally, the fact that the projective model structure is symmetric
monoidal follows from Corollary~2.21
of~\cite{hovey-christensen-relative}, the fact that $A$, the unit of
$\smash$, is a relative projective, and the fact that relative
projectives are closed under $\smash$.
\end{proof}

Note that, when $\hopfalg$ is discrete, a $\Gamma$-comodule is the same
thing as an $A$-module.  In this case, the relative projectives are just
the projective $A$-modules, and we see that the projective model
structure agrees with the usual projective model structure on $\Ch{A}$,
in which the fibrations are the surjections and the weak equivalences
are the homology isomorphisms.  

We point out that there is another model structure on $\Ch{\Gamma}$
given by~\cite[Example~3.4]{hovey-christensen-relative} called the
\emph{absolute model structure}.  In this model structure, the weak
equivalences are the chain homotopy equivalences, the cofibrations are
the degreewise split monomorphisms, and the fibrations are the
degreewise split epimorphisms.  Since the generating cofibrations of the
projective model structure are degreewise split monomorphisms, and the
generating trivial cofibrations are chain homotopy equivalences, we
conclude that the identity functor is a left Quillen functor from the
projective model structure to the absolute model structure.  In
particular, a trivial cofibration in the projective model structure is a
chain homotopy equivalence, and all chain homotopy equivalences are
projective equivalences.

The symmetric monoidal product behaves particularly well with respect to
the projective model structure.  

\begin{proposition}\label{prop-projective-monoid}
Suppose $\hopfalg$ is a flat Hopf algebroid.  Then the projective model
structure satisfies the monoid axiom.  Furthermore, if $X$ is 
cofibrant and $f$ is a projective equivalence, then $X\smash f$ is a
projective equivalence. 
\end{proposition}

This proposition is important because of the work of Schwede and
Shipley~\cite{schwede-shipley-monoids}, who introduced the monoid axiom.
As a consequence of their work and
Proposition~\ref{prop-projective-monoid}, given a monoid $R$ in
$\Ch{\Gamma }$, which is just a differential graded comodule algebra,
there is a model structure on (differential graded) $R$-modules in which
the fibrations are underlying projective fibrations and the weak
equivalences are underlying projective equivalences.  There is also a
similar model structure on differential graded comodule algebras, and a
projective equivalence $R\xrightarrow{}R'$ of differential graded
comodule algebras induces a Quillen equivalence from $R$-modules to
$R'$-modules.

\begin{proof}
The monoid axiom, introduced by Schwede and Shipley
in~\cite{schwede-shipley-monoids}, asserts that, if $K$ is the class of
maps $\{j\smash X \}$ where $j$ is a generating trivial cofibration and
$X$ is arbitrary, then all transfinite
compositions of pushouts of maps of $K$ are projective equivalences.
In the case at hand, $j$ is one of the maps $0\xrightarrow{}D^{n}P$,
where $P\in \cat{S}$.  It follows easily that $j\smash X$ is a
dimensionwise split monomorphism and a chain homotopy equivalence, so a
trivial cofibration in the absolute model structure.  Thus, all
transfinite compositions of pushouts of maps of $K$ are also trivial
cofibrations in the absolute model structure, and so in particular chain
homotopy equivalences.  Hence they are also projective equivalences.  

Now suppose $f\mathcolon Y\xrightarrow{}Z$ is a projective equivalence
and $X$ is cofibrant. We want to show that $X\smash f$ is a projective
equivalence.  Since $X$ is cofibrant, $X$ is a retract of a colimit of a
sequence of complexes $\{X_{i} \}_{i<\lambda }$, where
$X_{i}\xrightarrow{}X_{i+1}$ is a degreewise split monomorphism whose
cokernel $C_{i}$ is a complex of relative projectives with zero
differential.  Since projective equivalences are closed under filtered
colimits, it suffices to show that $X_{i}\smash f$ is a projective
equivalence for all $i\leq \lambda$, where $X_{\lambda}=\colim X_{i}$.
We prove this by transfinite induction on $i$.  We assume the base case
$i=0$ for the moment.  The limit ordinal case follows from the fact the
filtered colimits of projective equivalences are projective
equivalences.  For the sucessor ordinal case, we have the commutative
diagram below.
\[
\begin{CD}
0 @>>> X_{i} \smash Y @>>> X_{i+1}\smash Y @>>> C_{i} \smash Y @>>> 0 \\
@. @VX_{i}\smash fVV @VX_{i+1}\smash fVV @VVC_{i}\smash fV \\
0 @>>> X_{i}\smash Z @>>> X_{i+1}\smash Z @>>> C_{i}\smash Z @>>> 0
\end{CD}
\]
The rows of this diagram are degreewise split and short exact, since
$X_{i}\xrightarrow{}X_{i+1}$ is degreewise split.  Now apply the functor
$\Gamma \comod (P,-)$ to this diagram for a fixed $P\in \cat{S}$.  The
rows of the resulting diagram will still be short exact.  The induction
hypothesis tells us that $\Gamma \comod (P,X_{i}\smash f)$ is a homology
isomorphism, and the base case of the induction (which we have
postponed) tells us that $\Gamma \comod (P, C_{i}\smash f)$ is a
homology isomorphism.  The long exact sequence in homology then implies
that $\Gamma \comod (P, X_{i+1}\smash f)$ is a homology isomorphism.  

We are left with showing that $C\smash f$ is a projective equivalence,
where $C$ is a complex of relative projectives with zero differential.
Then $C\cong \bigoplus_{n} S^{n}C_{n}$.  Again using the fact that the
objects in $\cat{S}$ are finitely presented, we find that it suffices to
show that $S^{n}C_{n}\smash f$ is a projective equivalence.  But $C_{n}$
is a retract of a direct sum of elements of $\cat{S}$.  Another use of
the fact that objects in $S$ are finitely presented reduces us to
showing that $S^{n}Q\smash f$ is a projective equivalence, for $Q\in
\cat{S}$.  Assume $P\in \cat{S}$.  Then, since $Q$ is
strongly dualizable in $\Gamma \comod$ by
Proposition~\ref{prop-comodules-dualizable}, we have
\[
\Gamma \comod (P,S^{n}Q\smash f) \cong \Sigma^{-n}\Gamma \comod (P\smash DQ,f)
\]
which is a homology isomorphism since $P\smash DQ$ is also
(isomorphic to something) in $\cat{S}$.  
\end{proof}

An obvious drawback with the projective model structure is that is
difficult to tell what the weak equivalences look like.  We do have the
following proposition.  

\begin{proposition}\label{prop-projective-Adams}
Suppose $\hopfalg$ is a flat Hopf algebroid for which the dualizable
$\Gamma $-comodules generate the category of $\Gamma $-comodules.  Then
every projective fibration is surjective, and every projective
equivalence is a homology isomorphism.
\end{proposition}

\begin{proof}
Suppose $p\mathcolon X\xrightarrow{}Y$ is a projective fibration, and
$y\in Y_{n}$.  By Proposition~\ref{prop-generators}, there is a
comodule $P$ in $\cat{S}$ and a map $f\mathcolon P\xrightarrow{}Y_{n}$ whose
image contains $y$.  Suppose $f(t)=y$.  Since $p$ is a projective
fibration, there is a map $g\mathcolon P\xrightarrow{}X_{n}$ such that
$pg=f$.  In particular, $pg(t)=y$, so $p$ is surjective. 

Now suppose $p$ is a projective equivalence.  We wish to show that $p$
is a homology isomorphism.  Every projective trivial cofibration is a
chain homotopy equivalence, so a homology isomorphism.  We can thus
assume that $p$ is a projective trivial fibration.  In particular, $p$
is surjective.  Thus it suffices to show that $\ker p$ is exact.  We
know that $\ker p\xrightarrow{}0$ is a projective trivial fibration, so
$\Gamma \comod (P,\ker p)$ is exact for all $P\in \cat{S}$.  Suppose
that $x$ is a cycle in $\ker p_{n}$.  Let $Z_{n}$ denote the comodule of
cycles in $\ker p_{n}$.  Then there is a $P\in \cat{S}$, a $t\in P$, and
a comodule map $f\mathcolon P\xrightarrow{}Z_{n}$ such that $f(t)=x$, by
Proposition~\ref{prop-generators}.  The map $f$ is a cycle in $\Gamma
\comod (P,\ker p)$, so there is a map $g\mathcolon P\xrightarrow{}\ker
p_{n+1}$ such that $dg=f$.  In particular, $dg(t)=x$, so $x$ is a
boundary.
\end{proof}

\subsection{Naturality}\label{subsec-proj-naturality}

We now show that the projective model structure is natural in $\hopfalg
$.  

\begin{proposition}\label{prop-projective-natural}
Suppose $\Phi \mathcolon \hopfalg \xrightarrow{}\otherhopfalg$ is a map
of flat Hopf algebroids.  Then $\Phi$ induces a left Quillen functor
$\Phi_{*}\mathcolon \Ch{\Gamma}\xrightarrow{}\Ch{\Sigma}$ of the
projective model structures.  
\end{proposition}

\begin{proof}
We have seen in Proposition~\ref{prop-comodules-adjunction} that $\Phi$
induces an adjunction 
\[
(\Phi_{*},\Phi^{*})\mathcolon \Gamma \comod \xrightarrow{}\Sigma \comod.
\]
This prolongs to an adjunction
$(\Phi_{*},\Phi^{*})\mathcolon \Ch{\Gamma}\xrightarrow{}\Ch{\Sigma}$ by
defining $\Phi_{*}$ and $\Phi^{*}$ degreewise.  Since $\Phi _{*}$ is
symmetric monoidal, it preserves dualizable comodules.  This is easy to
see directly in this case, since if $M$ is finitely generated and
projective over $A$, then $B\otimes M$ is finitely generated and
projective over $B$.  It follows easily that $\Phi _{*}$ takes the
generating (trivial) cofibrations of the projective model structure on
$\Ch{\Gamma }$ to (trivial) cofibrations in the projective model
structure on $\Ch{\Sigma }$.
\end{proof}

Note that there is a map of Hopf algebroids $\Phi \mathcolon \hopfalg
\xrightarrow{}(A,A)$ which is the identity on $A$ and $\epsilon$ on
$\Gamma$.  The functor $\Phi_{*}$ is just the forgetful functor from
$\Gamma$-comodules to $A$-modules, and the right adjoint $\Phi^{*}$ is
the extended comodule functor.  Hence, if $f$ is a homology isomorphism
of complexes of $A$-modules, then $\Gamma \otimes f$ is a projective
equivalence.  

\begin{theorem}\label{thm-proj-weak}
Suppose $\Phi \mathcolon \hopfalg \xrightarrow{}\otherhopfalg$ is a weak
equivalence of flat Hopf algebroids.  Then $\Phi_{*}\mathcolon
\Ch{\Gamma}\xrightarrow{}\Ch{\Sigma}$ is a Quillen equivalence of the
projective model structures.  In fact, both $\Phi_{*}$ and $\Phi^{*}$
preserve and reflect projective equivalences.  
\end{theorem}

\begin{proof}
Since $\Phi ^{*}$ is a right Quillen functor, it preserves weak
equivalences (between fibrant objects, but everything is fibrant).
Since $\Phi _{*}$ is an equivalence, the unit $X\xrightarrow{}\Phi
^{*}\Phi _{*}X$ is an isomorphism, so $\Phi _{*}$ reflects projective
equivalences.  On the other hand, $\Phi^{*}$ is a symmetric monoidal
left adjoint since $\Phi_{*}$ is an equivalence of categories.  In
particular, $\Phi ^{*}$ preserves dualizable comodules.  Thus $\Phi
^{*}$ is a left (and right) Quillen functor of the projective model
structures.  Hence $\Phi_{*}$ is also a left and right Quillen functor,
so $\Phi _{*}$ preserves projective equivalences.  Thus $\Phi ^{*}$
reflects projective equivalences.  

To show that $\Phi _{*}$ is a Quillen equivalence, we need to show that,
for $X$ cofibrant and $Y$ fibrant, a map $f\mathcolon \Phi
_{*}X\xrightarrow{}Y$ is a projective equivalence if and only if its
adjoint $g\mathcolon X\xrightarrow{}\Phi ^{*}Y$ is a projective
equivalence.  Recall that $g$ is obtained from $f$ as the composite 
\[
X \xrightarrow{} \Phi ^{*}\Phi _{*}X \xrightarrow{\Phi ^{*}f} \Phi
^{*}Y. 
\]
The first map in this composite is an isomorphism.  Thus $g$ is a
projective equivalence if and only if $\Phi ^{*}f$ is so.  But $\Phi
^{*}$ preserves and reflects projective equivalences, so $\Phi ^{*}f$ is
a projective equivalence if and only if $f$ is so.  
\end{proof}

\subsection{The cobar resolution}\label{subsec-proj-cobar}

The projective model structure is clearly not the model structure we
want, because 
\[
\sho \Ch{\Gamma }(S^{0}A, S^{0}A)_{*} \cong A
\]
concentrated in degree $0$, because $A$ is both cofibrant and fibrant in
the projective model structure.  Recall that we want 
\[
\stable{\Gamma }(S^{0}A, S^{0}A)_{*} \cong \Ext _{\Gamma }^{*}(A,A).
\]
Therefore, we have to get an injective resolution of $A$ involved.

The injective resolution we choose is the cobar
resolution~\cite[A1.2.11]{ravenel}, though we offer a simpler
construction of it.  Suppose $M$ is a $\Gamma$-comodule.  Then
$\psi$ is a natural comodule embedding $M\xrightarrow{}\Gamma \otimes M$
of $M$ into an extended comodule, which is split over $A$ by $\epsilon
\otimes 1$.  We can iterate this to construct a resolution of $M$ by
extended $A$-comodules.  The most important case is when $M=A$.  We
begin with the $A$-split short exact sequence of comodules 
\[
0\xrightarrow{} A\xrightarrow{\eta_{L}}\Gamma \xrightarrow{}
\overline{\Gamma} \xrightarrow{}0.
\]
Here $\overline{\Gamma}$ is of course the cokernel of $\eta_{L}$, but it
is easily seen to be isomorphic to $\ker \epsilon$.  When we think of it
as $\ker \epsilon$, the coaction is defined by $\psi (x)=\Delta
(x)-x\otimes 1$.  We can then tensor this sequence with
$\overline{\Gamma}^{\smash s}$ to get the $A$-split short exact sequence
of comodules 
\[
0\xrightarrow{} \overline{\Gamma}^{\smash s} \xrightarrow{}\Gamma \smash
\overline{\Gamma}^{\smash s} \xrightarrow{} \overline{\Gamma}^{\smash
(s+1)} \xrightarrow{} 0.  
\]
We splice these short exact sequences together to obtain a complex $LA$,
where $(LA)_{-n}=\Gamma \smash \overline{\Gamma}^{\smash n}$ for $n\geq
0$ and $(LA)_{-n}=0$ for $n<0$, and the differential is the composite
\[
\Gamma \smash \overline{\Gamma}^{\smash n} \xrightarrow{}
\overline{\Gamma}^{\smash (n+1)} \xrightarrow{} \Gamma \smash
\overline{\Gamma}^{\smash (n+1)}.  
\]
In particular, there is a homology isomorphism $S^{0}A\xrightarrow{}LA$
induced by $\eta_{L}$, so that $LA$ is a resolution of $A$, and the
cycle comodule $Z_{-n}(LA)$ is isomorphic to $\overline{\Gamma}^{\smash
n}$ for $n>0$ (and $A$ for $n=0$).  Furthermore, the $A$-splittings
patch together to show that $S^{0}A\xrightarrow{}LA$ is a chain homotopy
equivalence of complexes of $A$-modules.  

The complex $LA$ will be very important in the rest of this
paper, but $LA$ is not cofibrant in the projective model structure, since
$(LA)_{0}=\Gamma$ is not even projective over $A$ in general.  The
following proposition is then crucial for us.  

\begin{proposition}\label{prop-projective-LA}
Suppose that $\hopfalg$ is a flat Hopf algebroid, that the dualizable
$\Gamma $-comodules generate the category of $\Gamma $-comodules, and
that $\overline{\Gamma }\smash X$ is projectively trivial when $X$ is
so.  Let $LA$ denote the cobar resolution of $A$.
\begin{enumerate}
\item [(a)] If $p$ is a projective fibration, then $LA\smash p$ is a
projective fibration.  
\item [(b)] If $p$ is a projective equivalence, then $LA\smash p$ is a
projective equivalence.
\end{enumerate}
\end{proposition}

Note that, for part~(a), it is sufficient to assume that dualizable
comodules generate.

In view of Proposition~\ref{prop-projective-LA}, we make the following
definition.  

\begin{definition}\label{defn-amenable}
A Hopf algebroid $\hopfalg $ is \emph{amenable} when it is flat, the
dualizable $\Gamma $-comodules generate the category of $\Gamma
$-comodules, and $\overline{\Gamma }\smash (-)$ preserves projectively
trivial complexes.  
\end{definition}

For Proposition~\ref{prop-projective-LA} to be of use, we need to know
that amenable Hopf algebroids do exist. 

\begin{proposition}\label{prop-amenable}
Every Adams Hopf algebroid is amenable.  
\end{proposition}

The rest of this section will be devoted to proving
Propositions~\ref{prop-projective-LA} and~\ref{prop-amenable}.
Proposition~\ref{prop-amenable} is an immediate consequence of the
following lemma and Proposition~\ref{prop-amenable-Adams}.  

\begin{lemma}\label{lem-proj-amenable}
Suppose $\hopfalg $ is a flat Hopf algebroid. 
\begin{enumerate}
\item [(a)] If $M$ is a filtered colimit of dualizable comodules, then
$M\smash (-)$ preserves projectively trivial complexes. 
\item [(b)] If $\hopfalg $ is Adams, then $\overline{\Gamma }$ is a
filtered colimit of dualizable comodules.  
\end{enumerate}
\end{lemma}

\begin{proof}
For part~(a), recall that filtered colimits of projective equivalences
are projective equivalences.  We can therefore assume that $M$ itself is
a dualizable comodule.  Suppose then that $X$ is projectively trivial.
We must show that $\Gamma \comod (P, M\smash X)$ has no homology for all
dualizable comodules $P$.  But 
\[
\Gamma \comod (P, M\smash X) \cong \Gamma \comod (P\smash DM,X)
\] 
since $M$ is dualizable.  Furthermore, $P\smash DM$ is again dualizable,
so since $X$ is projectively trivial, we are done.  

For part~(b), since $\hopfalg $ is an Adams Hopf algebroid, we have
$\Gamma =\colim_{i\in \cat{I}} \Gamma_{i}$ for a filtered small category
$\cat{I}$ of arrows $i\mathcolon \Gamma_{i}\xrightarrow{}\Gamma $ such
that $\Gamma_{i}$ is dualizable.  Let $\cat{J}$ denote the category of
factorizations $A\xrightarrow{}\Gamma_{i}\xrightarrow{i}\Gamma$ of
$\eta_{L}$ through an arrow of $\cat{I}$.  By abuse of notation, we
write the map $A\xrightarrow{}\Gamma_{i}$ as $\eta_{L}$ as well; note
that this $\eta_{L}$ must be a split monomorphism of $A$-modules, since
the usual $\eta_{L}$ is so.  We claim that $\cat{J}$ is filtered and
that the obvious functor from $\cat{J}$ to $\cat{I}$ is cofinal (see
Definition~2.3.8 of~\cite{hovey-axiomatic} for a reminder of what this
means).  This is a straightforwad consequence of the fact that $A$ is
itself finitely presented as a $\Gamma$-comodule.  It follows then that
$\colim_{j\in \cat{J}}\Gamma_{j}\cong \Gamma$, and therefore that
$\colim \overline{\Gamma_{j}}\cong \overline{\Gamma}$, where
$\overline{\Gamma_{j} }=\coker \eta_{L}$.  Each $\overline{\Gamma _{j}}$
is finitely generated and projective over $A$, and hence dualizable.
\end{proof}

Part~(a) of Proposition~\ref{prop-projective-LA} is an immediate
consequence of the following lemma and
Proposition~\ref{prop-projective-Adams}. 

\begin{lemma}\label{lem-projective-LA-fib}
Let $\hopfalg $ be a flat Hopf algebroid. 
\begin{enumerate}
\item [(a)] If $M$ is an extended comodule, then $M\smash (-)$ takes
surjections of comodules to relative epimorphisms.  
\item [(b)] If $X$ is a complex of extended comodules, then $X\smash
(-)$ takes surjections of complexes to projective fibrations. 
\end{enumerate}
\end{lemma}

\begin{proof}
Suppose $f$ is a surjection of complexes, and $M\cong \Gamma \otimes N$
is an extended comodule.  Let $P$ be a dualizable comodule.  Then, using
Lemma~\ref{lem-different-tensor}, we
have
\begin{gather*}
\Gamma \comod (P, M\smash f) \cong \Gamma \comod (P, (\Gamma \otimes N)
\smash f) \\ 
\cong \Gamma \comod (P, \Gamma \otimes (N\otimes f)) \cong \Hom
_{A}(P,N\otimes f).
\end{gather*}
Since $f$ is surjective, so is $N\otimes f$.  Since $P$ is projective
over $A$, $\Hom _{A}(P, N\otimes f)$ is also surjective, so $M\smash f$
is a relative epimorphism.

Now suppose $X$ is a complex of extended comodules.  Then, in degree
$n$, we have 
\[
(X\smash f)_{n} \cong \bigoplus_{m} X_{m} \smash f_{n-m}. 
\]
Each map $f_{n-m}$ is surjective, so part~(a) assures us that
$X_{m}\smash f_{n-m}$ is a relative epimorphism.  Since the dualizable
complexes are finitely presented, direct sums of relative epimorphisms
are again relative epimorphisms.  Hence $X\smash f$ is a projective
fibration.  
\end{proof}

We are left with proving part~(b) of
Proposition~\ref{prop-projective-LA}. Our approach is similar to that of
Lemma~\ref{lem-projective-LA-fib}.  

\begin{lemma}\label{lem-projective-LA-equiv}
Suppose $\hopfalg $ is a flat Hopf algebroid. 
\begin{enumerate}
\item [(a)] Suppose $N$ is a flat $A$-module.  Then $(\Gamma \otimes
N)\smash (-)$ takes exact complexes to projectively trivial complexes. 
\item [(b)] Suppose $X$ is a bounded below complex such that
$X_{n}\smash (-)$ preserves projectively trivial complexes for all $n$.
Then $X\smash (-)$ preserves projectively trivial complexes. 
\item [(c)] Suppose $X$ is a complex of comodules such that $X_{n}\smash
(-)$ and $Z_{n}X\smash (-)$ preserve projectively trivial complexes for
all $n$.  Then $X\smash (-)$ preserves projectively trivial complexes.   
\end{enumerate}
\end{lemma}

Here $Z_{n}X$ denotes the cycles in degree $n$, as usual.  

\begin{proof}
For part~(a), suppose $Y$ is a projectively trivial complex and $P$ is a
dualizable comodule.  Then, using Lemma~\ref{lem-different-tensor}, we
have 
\[
\Gamma \comod (P, (\Gamma \otimes N)\smash Y) \cong \Gamma \comod (P,
\Gamma \otimes (N\otimes Y)) \cong \Hom _{A}(P, N\otimes Y).  
\]
Since $Y$ is exact and $N$ is flat, $N\otimes Y$ is also exact.  Since
$P$ is projective over $A$, $\Hom _{A}(P, N\otimes Y)$ is also exact, as
required.  

For part~(b), let $Y$ be a projectively trivial complex.  Without loss of
generality, we can assume that $X_{n}=0$ for $n<0$.  Suppose $P\in
\cat{S}$, and $z\mathcolon P\xrightarrow{}(X\smash Y)_{n}$ is a cycle in
$\Gamma \comod (P,X\smash Y)$.  Since $P$ is finitely presented as a
$\Gamma$-comodule,
\[
\Gamma \comod (P,(X\smash Y)_{n})\cong \bigoplus_{i=0}^{\infty } \Gamma
\comod (P,X_{i}\smash Y_{n-i}).
\]
We can therefore write $z=(z_{0},z_{1},\dotsc ,z_{i},\dotsc)$, where
$z_{i}\mathcolon P\xrightarrow{}X_{i}\smash Y_{n-i}$ and $z_{i}=0$ for
large $i$. Define the \emph{degree} of $z$ to be the largest $i$ such
that $z_{i}$ is nonzero.  We will show that every cycle $z$ is
homologous to a cycle of smaller degree; since there are no cycles of
degree $-1$ this will complete the proof.  Indeed, suppose $z$ has
degree $k$.  Then $z_{k}$ has to be a cycle in the complex $X_{k}\smash
Y$.  By assumption, $X_{k}\smash Y$ is projectively trivial, so $z_{k}$
must be a boundary in this complex.  This means that there is a
$w\mathcolon P\xrightarrow{}X_{k}\smash Y_{n-k+1}$ such that $(1\smash
d)w=z_{k}$.  But then $z$ is homologous to $z'=z+(-1)^{k+1}dw$, and one
can easily check that $w$ has degree $<k$.

For part~(c), again assume that $Y$ is projectively trivial.  Let
$X^{i}$ be the subcomplex of $X$ such that $X^{i}_{n}=X_{n}$ for $n>-i$,
$X^{i}_{n}=0$ for $n<-i$, and $X^{i}_{-i}=Z_{-i}X$.  By part~(b), each
of the complexes $X^{i}\smash Y$ is projectively trivial.  But $X=\colim
X^{i}$, so $X\smash Y=\colim X^{i}\smash Y$.  Since filtered colimits of
projective equivalences are projective equivalences, $X\smash Y$ is
therefore projectively trivial.  
\end{proof}

We can now prove Proposition~\ref{prop-projective-LA}(b). 

\begin{proof}[Proof of Proposition~\ref{prop-projective-LA}\ulp b\urp ]

We need to show that $LA \smash (-)$ preserves projective equivalences.
Since the projective trivial cofibrations are in particular chain
homotopy equivalences, $LA\smash (-)$ certainly takes them to projective
equivalences.  It therefore suffices to show that $LA\smash (-)$
preserves projective trivial fibrations.  By part~(a), $LA\smash (-)$
preserves projective fibrations, so it suffices to show that $LA\smash
(-)$ preserves projectively trivial complexes.  

In view of Lemma~\ref{lem-projective-LA-equiv}, it suffices to show that
$(LA)_{n}\smash (-)$ and $Z_{n}LA\smash (-)$ preserve projectively
trivial complexes.  Since $(LA)_{-n}=\Gamma \smash \overline{\Gamma
}^{\smash n}$ for $n\geq 0$, and $\overline{\Gamma }$ is flat as an
$A$-module since $\Gamma $ is so, Lemma~\ref{lem-projective-LA-equiv}(a)
guarantees that $(LA)_{n}\smash (-)$ preserves projectively trivial
complexes.  On the other hand, $Z_{-n}(LA)\cong \overline{\Gamma
}^{\smash n}$ for $n\geq 0$.  The amenable assumption guarantees that
$\overline{\Gamma }\smash (-)$ preserves projectively trivial complexes,
and then iteration shows that $Z_{-n}(LA)\smash (-)$ does so as well.
\end{proof}

\section{Homotopy groups}\label{sec-groups}

When the dualizable $\Gamma $-comodules generate the category of $\Gamma
$-comodules, we know from Proposition~\ref{prop-projective-Adams} that
projective equivalences are homology isomorphisms.  But
homology is not the most important functor of complexes of comodules;
homotopy is.  In this section we define and study the homotopy groups of
a chain complex of comodules.  We show that these homotopy groups are
closely related to $\Ext $ in the category of $\Gamma $-comodules and
have similar properties.  When $\hopfalg $ is amenable, every projective
equivalence is a homotopy isomorphism and every homotopy isomorphism is
a homology isomorphism.  The object of Section~\ref{sec-homotopy} will
then be to construct a model structure on $\Ch{\Gamma}$ in which the
weak equivalences are the homotopy isomorphisms.

\subsection{Relatively injective comodules}\label{subsec-injective}

To explain homotopy groups, we need to remind the reader of some of the
basic results on relatively injective comodules.  Some of this
can be found in~\cite[Appendix~1]{ravenel}.

\begin{definition}\label{defn-rel-inj}
Suppose $\hopfalg$ is a flat Hopf algebroid.  A comodule $I$ is called
\emph{relatively injective} if $\Gamma \comod (-,I)$ takes $A$-split
short exact sequences to short exact sequences.  
\end{definition}

\begin{lemma}\label{lem-rel-inj}
Suppose $\hopfalg$ is a flat Hopf algebroid.  The relatively injective
comodules are the retracts of extended comodules.  In particular, there
is a natural $A$-split embedding of any comodule into a relatively
injective comodule.  
\end{lemma}

\begin{proof}
We have $\Gamma \comod (-, \Gamma \otimes N)\cong \Hom _{A}(-,N)$.  Thus
extended comodules, and so also retracts of extended comodules, are
relatively injective.  Conversely, if $I$ is a relative injective, the
map $I\xrightarrow{\psi}\Gamma \otimes I$ must have a retraction, since
it is a map of comodules that is split over $A$ by $\epsilon \otimes 1$.
Thus $I$ is a retract of $\Gamma \otimes I$.  The natural $A$-split
embedding of the statement of the lemma is just $\psi \mathcolon
M\xrightarrow{}\Gamma \otimes M$.
\end{proof}

Of course, there are (absolutely) injective comodules as well.  A
similar argument shows that the injective comodules are retracts of
extended comodules $\Gamma \otimes I$, where $I$ is an injective
$A$-module.  But relatively injective comodules are much easier to work
with than injective comodules, partly because injective $A$-modules are
complicated, and partly because of the following lemma. 

\begin{lemma}\label{lem-rel-inj-closure}
Suppose $\hopfalg$ is a flat Hopf algebroid.
\begin{enumerate}
\item [(a)] Relatively injective comodules are closed under coproducts
and products.  
\item [(b)] If $M$ is an arbitrary comodule and $I$ is relatively
injective, then $I\smash M$ and $F(M,I)$ are relatively injective.  
\end{enumerate}
\end{lemma}

\begin{proof}
For part~(a), it suffices to show that extended comodules are closed
under coproducts and products.  But we have 
\[
\bigoplus (\Gamma \otimes M_{i})\cong \Gamma \otimes (\bigoplus M_{i})
\text{ and} \prod^{\Gamma } (\Gamma \otimes M_{i}) \cong \Gamma \otimes
(\prod M_{i}),
\]
the latter by the construction of products in
Proposition~\ref{prop-comod-products}.  

For part~(b), we first prove that $F(M,I)$ is relatively injective.  We
must show that $\Gamma \comod (-, F(M,I))$ takes $A$-split short exact
sequences to short exact sequences.  But $\Gamma \comod (-,F(M,I))$ is
naturally isomorphic to $\Gamma \comod (-\smash M,I)$, so this is clear.

To show that $I\smash M$ is relatively injective, note that $I\smash M$
is a retract of $(\Gamma \otimes I)\smash M$, which is isomorphic to
$\Gamma \smash I\smash M$ by Lemma~\ref{lem-different-tensor}.  On the
other hand, another use of Lemma~\ref{lem-different-tensor} shows that
$\Gamma \smash I\smash M$ is isomorphic to the extended comodule $\Gamma
\otimes (I\otimes M)$, completing the proof.  
\end{proof}

Relatively injective comodules can be used to compute $\Ext$ when the
source is projective over $A$. 

\begin{lemma}\label{lem-proj-inj}
Suppose $\hopfalg$ is a flat Hopf algebroid, $P$ is a $\Gamma$-comodule
that is projective over $A$, and $I$ is a relatively injective
comodule.  Then $\Ext^{n}_{\Gamma}(P,I)=0$ for all $n>0$.  Hence, if
$I_{*}$ is a resolution of $M$ by relatively injective comodules, 
\[
\Ext^{n}_{\Gamma}(P,M) \cong H_{-n}(\Gamma \comod (P,I_{*})).
\]
\end{lemma}

\begin{proof}
The second statement is an immediate consequence of the first.  It
suffices to prove the first statement for $I=\Gamma \otimes N$, since
every relatively injective comodule is a retract of an extended
comodule.  Let $J_{*}$ be an injective resolution of $N$ in the category
of $A$-modules.  Then $\Gamma \otimes J_{*}$ is an injective resolution
of $\Gamma \otimes N$ in $\Gamma \comod$, since $\Gamma$ is flat over
$A$.  Hence
\[
\Ext^{n}_{\Gamma}(P,\Gamma \otimes N)\cong H_{-n}(\Gamma \comod (P,
\Gamma \otimes J_{*})) \cong H_{-n}\Hom_{A}(P,J_{*}) \cong
\Ext^{n}_{A}(P,N).
\]
Since $P$ is projective over $A$, this group is $0$ for $n>0$.  
\end{proof}

\subsection{Homotopy groups}\label{subsec-groups-homotopy}

Now recall that $LA$ denotes a specific resolution of $A$ by relatively
injective comodules, defined in Section~\ref{subsec-proj-cobar}, such
that the map $S^{0}A\xrightarrow{}LA$ is a chain homotopy equivalence
over $A$.  It follows that $LA\smash M$ is a resolution of $M$ for any
comodule $M$.  It is in fact a resolution by relative injectives by
Lemma~\ref{lem-rel-inj-closure}.  Thus we have
\[
\Ext^{n}_{\Gamma}(P,M) \cong H_{-n}(\Gamma \comod (P,LA\smash M))
\]
for any comodule $P$ that is projective over $A$.  

We now extend this definition, replacing $M$ by a complex $X$.  

\begin{definition}\label{defn-homotopy}
Suppose $\hopfalg$ is a flat Hopf algebroid, $X\in \Ch{\Gamma}$, $P\in
\cat{S}$, and $n\in \Z$.  Define the
\emph{$n$th homotopy group of $X$ with coefficients in $P$},
$\pi_{n}^{P}(X)$, by $\pi_{n}^{P}(X)=H_{-n}(\Gamma \comod (P,LA\smash X))$.  
\end{definition}

We need to say a few words about grading.  We have essentially two
choices; we can grade homotopy as if it were the homotopy groups of a
space, or we can grade it as if it were the $\Ext $ groups of a
comodule.  Either way has problems; grading it like $\Ext $ means the
exact sequences on homotopy go up instead of down in dimension, but
grading it like homotopy means the homotopy groups of $A$ will be
concentrated in negative degrees.  Following Palmieri's work on the
Steenrod algebra~\cite{palmieri-book}, we choose to grade it like $\Ext
$.  This extends to bigrading as well; if $\hopfalg $ is a graded Hopf
algebroid, as it always is in algebraic topology, we define
\[
\pi _{s,t}^{P}(X) = H_{-s,t}(\Gamma \comod (P, LA\smash X)).  
\]

The homotopy groups are of course functorial in $X$, and they satisfy
the expected properties, correcting for the strange grading. 

\begin{lemma}\label{lem-homotopy-properties}
Suppose $\hopfalg$ is a flat Hopf algebroid.
\begin{enumerate}
\item [(a)] A short exact sequence of complexes 
\[
0 \xrightarrow{} X \xrightarrow{} Y \xrightarrow{} Z \xrightarrow{} 0
\]
induces a natural long exact sequence 
\[
\dotsb \xrightarrow{} \pi_{n}^{P}(X) \xrightarrow{} \pi_{n}^{P}(Y)
\xrightarrow{} \pi_{n}^{P}(Z) \xrightarrow{} \pi_{n+1}^{P}(X)
\xrightarrow{} \dotsb
\]
\item [(b)] If $X$ is a filtered colimit of complexes $X^{i}$, then
$\pi_{n}^{P}(X)\cong \colim \pi_{n}^{P}(X^{i})$.  
\end{enumerate}
\end{lemma}

\begin{proof}
For part~(a), since $LA$ is degreewise flat over $A$, the sequence
\[
0 \xrightarrow{} LA\smash X \xrightarrow{} LA\smash Y\xrightarrow{}
LA\smash Z \xrightarrow{} 0
\]
remains exact.  By Lemma~\ref{lem-rel-inj-closure}, $LA\smash X$ is a
complex of relative injectives.  Therefore, the sequence of complexes
\[
0 \xrightarrow{} \Gamma \comod (P, LA\smash X) \xrightarrow{} \Gamma
\comod (P, LA\smash Y) \xrightarrow{} \Gamma \comod (P,LA\smash Z)
\xrightarrow{} 0
\]
remains exact, by Lemma~\ref{lem-proj-inj}.  The long exact sequence in
homology of this short exact sequence finishes the proof of part~(a).  

For part~(b), we simply note that $LA\smash -$, $\Gamma \comod (P,-)$,
and homology all commute with filtered colimits.  
\end{proof}

\subsection{Homotopy isomorphisms}\label{subsec-groups-iso}

A chain map $\phi$ is called a \emph{homotopy isomorphism} if
$\pi_{n}^{P}(\phi)$ is an isomorphism for all $n\in \Z$ and all $P\in
\cat{S}$.  Note that $\phi $ is a homotopy isomorphism if and only if
$LA\smash \phi$ is a projective equivalence.  We claim that homotopy
isomorphisms are the natural notion of weak equivalence in
$\Ch{\Gamma}$.

\begin{proposition}\label{prop-homotopy-projective}
Suppose $\hopfalg$ is an amenable Hopf algebroid.  Then every projective
equivalence is a homotopy isomorphism, and every homotopy isomorphism is
a homology isomorphism.
\end{proposition}
  
\begin{proof}
Suppose $p$ is a projective equivalence.  Then
Proposition~\ref{prop-projective-LA} tells us that $LA\smash p$ is also
a projective equivalence, so $p$ is a homotopy isomorphism.  Now suppose
$p\mathcolon X\xrightarrow{}Y$ is a homotopy isomorphism.  Then
$LA\smash p$ is a projective equivalence, and hence a homology
isomorphism by Proposition~\ref{prop-projective-Adams}.  But
$A\xrightarrow{}LA$ is a chain homotopy equivalence over $A$, so
$X\xrightarrow{}LA\smash X$ and $Y\xrightarrow{}LA\smash Y$ are also
chain homotopy equivalences over $A$, and in particular homology
isomorphisms.  Hence $p$ is a homology isomorphism.
\end{proof}

Homotopy isomorphisms have the properties one would hope for in a
collection of weak equivalences.

\begin{proposition}\label{prop-homotopy-properties}
Suppose $\hopfalg$ is a flat Hopf algebroid.
\begin{enumerate}
\item [(a)] Homotopy isomorphisms are closed under retracts and have the
two out of three property.  
\item [(b)] Homotopy isomorphisms are closed under filtered colimits.  
\item [(c)] If $f$ is an injective homotopy isomorphism,
and $g$ is a pushout of $f$, then $g$ is an injective homotopy
isomorphism.  Dually, if $f$ is surjective homotopy isomorphism, and $g$
is a pullback of $f$, then $g$ is a surjective homotopy isomorphism.  
\item [(d)] If $f$ is a homotopy isomorphism, then any pushout of $f$
through an injective map is again a homotopy isomorphism.  Dually, any
pullback of $f$ through a surjective map is again a homotopy
isomorphism.  
\item [(e)] Suppose $f\mathcolon X\xrightarrow{}Y$ is an injective
homotopy isomorphism, and $g\mathcolon A\xrightarrow{}B$ is a
cofibration.  Then the pushout product
\[
f\boxprod g\mathcolon (X\smash B) \amalg_{X\smash A} (Y\smash A)
\xrightarrow{} Y\smash B
\]
is an injective homotopy isomorphism.  
\end{enumerate}
\end{proposition}

\begin{proof}
We leave part~(a) to the reader.  Part~(b) is immediate from the fact
that homotopy groups commute with filtered colimits.  For part~(c),
suppose $g$ is a pushout of the injective homotopy isomorphism $f$.
Then $g$ is injective, with cokernel $\coker f$.  Since $f$ is a
homotopy isomorphism, the long exact sequence of
Lemma~\ref{lem-homotopy-properties} shows that $\coker f$ has zero
homotopy.  Another use of that long exact sequence shows that $g$ is a
homotopy isomorphism.  The dual case is similar.  

For part~(d), suppose that $g\mathcolon B\xrightarrow{}D$ is the pushout
of the homotopy isomorphism $f\mathcolon A\xrightarrow{}C$ through the
injection $i\mathcolon A\xrightarrow{}B$.  Then we have the map of short
exact sequences below.
\[
\begin{CD}
0 @>>> A @>i>> B @>>> X @>>> 0 \\
@. @VfVV @VgVV @| \\
0 @>>> C @>>> D @>>> X @>>> 0
\end{CD}
\]
The long exact sequence in homotopy and the five lemma show that $g$ is
a homotopy isomorphism.  The dual case is similar.  

For part~(e), note that parts~(b) and~(c) imply that injective homotopy
isomorphisms are closed under pushouts and filtered colimits, hence
transfinite compositions.  Thus it suffices to prove part~(e) when $g$
is one of the generating cofibrations $S^{n-1}P\xrightarrow{}D^{n}P$ of
the projective model structure, by Lemma~4.2.4 of~\cite{hovey-model}.
We leave to the reader the check that $f\boxprod g$ is injective in this
case, and just prove it is a homotopy isomorphism.  Since $f$ is a
homotopy isomorphism, $LA\smash f$ is a projective equivalence.
Therefore, $(LA\smash f)\smash S^{n-1}P$ is a projective equivalence by
Proposition~\ref{prop-projective-monoid}, and so $f\smash S^{n-1}P$ is a
homotopy isomorphism.  Similarly, $f\smash D^{n}P$ is a homotopy
isomorphism.  Both such maps are also injective, since $P$ is flat over
$A$.  Part~(c) implies that the pushout
\[
X \smash D^{n}P \xrightarrow{} (X\smash D^{n}P)\amalg_{X\smash S^{n-1}P}
(Y\smash S^{n-1}P)
\]
is also an injective homotopy isomorphism.  The two out of three
property for homotopy isomorphisms then implies that $f\boxprod g$ is a
homotopy isomorphism.  
\end{proof}

Our next goal is to give some useful examples of homotopy isomorphisms
that are not projective equivalences.  We begin with the following
lemma.

\begin{lemma}\label{lem-contractible-to-injective}
Let $\hopfalg $ be a flat Hopf algebroid.  Suppose $X\in \Ch{\Gamma }$
is bounded above and contractible as a complex of $A$-modules, and $Y\in
\Ch{\Gamma }$ is a complex of relatively injective comodules.  Then
every chain map $f\mathcolon X\xrightarrow{}Y$ is chain homotopic to
$0$.
\end{lemma}

\begin{proof}
We construct a chain homotopy $D_{n}\mathcolon
X_{n}\xrightarrow{}Y_{n+1}$ by downwards induction on $n$.  Getting
started is easy, since $X_{n}=0$ for large $n$.  Suppose we have
constructed $D_{n+1}$ and $D_{n+2}$ such that
$dD_{n+2}+D_{n+1}d=f_{n+2}$.  We need to construct $D_{n}$ such that
$dD_{n+1}+D_{n}d=f_{n+1}$.  One can readily verify that 
\[
(f_{n+1}-dD_{n+1})d =0
\]
and so $f_{n+1}-dD_{n+1}$ defines a map $g\mathcolon X_{n+1}/\im d =
X_{n+1}/\ker d \xrightarrow{}Y_{n+1}$.  On the other hand, we are given
that $X$ is $A$-contractible, so there are $A$-module maps
$s_{n}\mathcolon X_{n}\xrightarrow{}X_{n+1}$ such that $ds+sd=1$.  In
particular, $d\mathcolon X_{n+1}/\ker d\xrightarrow{}X_{n}$ is an
$A$-split monomorphism.  Since $Y_{n+1}$ is relatively injective, there
is a map $D_{n}\mathcolon X_{n}\xrightarrow{}Y_{n+1}$ such that
$D_{n}d=f_{n+1}-dD_{n+1}$.  This completes the induction step and the
proof.  
\end{proof}

This gives the following proposition.

\begin{proposition}\label{prop-homotopy-none}
Let $\hopfalg$ be a flat Hopf algebroid, and suppose $f\mathcolon
X\xrightarrow{}Y$ is a map of bounded above complexes in $\Ch{\Gamma}$
that is an $A$-split monomorphism in each degree and a chain homotopy
equivalence of complexes of $A$-modules.  Then $LA\smash f$ is a chain
homotopy equivalence.  In particular, $f$ is a homotopy isomorphism.
\end{proposition}

\begin{proof}
Let $Z$ denote the cokernel of $f$.  Then $Z$ is bounded above and
contractible as a complex of $A$-modules (one can check this directly,
but it also follows because $f$ is a trivial cofibration in the absolute
model structure on
$\Ch{A}$~\cite[Example~3.4]{hovey-christensen-relative}).  Therefore
$LA\smash Z$ is a bounded above complex of relatively injective
comodules that is contractible over $A$.
Lemma~\ref{lem-contractible-to-injective} implies that $LA\smash Z$ is
contractible.  Since $LA$ is degreewise flat over $A$, $LA\smash Z$ is
the cokernel of $LA\smash f$.  Furthermore, $LA\smash f$ is a degreewise
$A$-split monomorphism of relatively injective comodules, so it is a
degreewise split monomorphism.  It follows that $LA\smash f$ is a chain
homotopy equivalence.  
\end{proof}

\begin{corollary}\label{cor-homotopy-none}
Suppose $\hopfalg$ is a flat Hopf algebroid.  Then the map
\[
\eta_{L}\smash X\mathcolon X\xrightarrow{} LA\smash X
\]
is a homotopy isomorphism for all complexes $X$.  In
fact, $LA\smash \eta_{L}\smash X$ is a chain homotopy equivalence.  
\end{corollary}

\begin{proof}
The map $\eta_{L}\mathcolon A\xrightarrow{}LA$ is a map of bounded above
complexes that is a degreewise $A$-split monomorphism and an $A$-chain
homotopy equivalence.  Proposition~\ref{prop-homotopy-none} implies that
$LA\smash \eta_{L}$ is a chain homotopy equivalence.  One can easily
check that this forces $LA\smash \eta_{L}\smash X$ to be a chain
homotopy equivalence for any $X$, and so $\eta_{L}\smash X$ is a
homotopy isomorphism.
\end{proof}

\section{Localization}\label{sec-model}

In the next section, we will localize the projective model structure to
obtain a model structure on $\Ch{\Gamma}$ in which the weak equivalences
are the homotopy isomorphisms.  To prove that this construction works,
we need some general results about Bousfield localization of model
categories.  The basic reference for Bousfield localization
is~\cite{hirschhorn}, but the results we prove in this section are
new. 

Suppose we have a model category $\cat{M}$ and a class of maps
$\cat{T}$.  The Bousfield localization $L_{\cat{T}}\cat{M}$ of $\cat{M}$
with respect to $\cat{T}$ is a new model structure on $\cat{M}$, with the
same cofibrations as the given one, in which the maps of $\cat{T}$ are
weak equivalences.  Futhermore, it is the initial such model category,
in the sense that if $F\mathcolon \cat{M}\xrightarrow{}\cat{N}$ is a
(left) Quillen functor that sends the maps of $\cat{T}$ to weak
equivalences, then $F\mathcolon L_{\cat{T}}\cat{M}\xrightarrow{}\cat{N}$
is also a Quillen functor.  

The Bousfield localization is known to exist when $\cat{T}$ is a set,
$\cat{M}$ is left proper, and, in addition, $\cat{M}$ is either
cellular~\cite{hirschhorn} or combinatorial (unpublished work of Jeff
Smith).  The cellular condition is technical, but has the virtue of
being written down and of applying to topological spaces.  The
combinatorial condition is simpler; it just means that $\cat{M}$ is
cofibrantly generated and locally presentable as a category.

To describe the localized model structure, it is necessary to recall
that any model category $\cat{M}$ possesses a unital action by the
category $\SSet$ of simplicial sets.  That is, there is a bifunctor
$\cat{M}\times \SSet \xrightarrow{}\cat{M}$ that takes $(X,K)$ to
$X\otimes K$ described in~\cite[Chapter~5]{hovey-model}.  This is a
unital action but is not associative; it induces an associative action
of $\sho \SSet$ on $\sho \cat{M}$.  In fact, $\sho \cat{M}$ is not only
tensored over $\sho \SSet$, but also cotensored and enriched over $\sho
\SSet $~\cite[Chapter~5]{hovey-model}.  The enrichment is denoted by
$\map (X,Y)\in \SSet$.  

Now, a fibrant object $X$ in $L_{\cat{T}}\cat{M}$, called a
\emph{$\cat{T}$-local fibrant object}, is a fibrant object $X$ in
$\cat{M}$ such that $\map (f,X)$ is a weak equivalence of simplicial
sets for all $f\in \cat{T}$.  Adjointness gives an equivalent
description, as follows.  Given a map $f\mathcolon X\xrightarrow{}Y$ in
$\cat{M}$, let $\widetilde{f}\mathcolon
\widetilde{X}\xrightarrow{}\widetilde{Y}$ denote a cofibration that is a
cofibrant approximation to $f$.  This means that $\widetilde{X}$ and
$\widetilde{Y}$ are cofibrant, $\widetilde{f}$ is a cofibration, and we
have the commutative diagram below
\[
\begin{CD}
\widetilde{X} @>\widetilde{f}>> \widetilde{Y} \\
@VVV @VVV \\
X @>>f> Y
\end{CD}
\]
where the vertical arrows are weak equivalences.  Then a \emph{horn on
$f$} is one of the maps $\widetilde{f}\boxprod i_{n}$ for $n\geq 0$,
where $i_{n}\mathcolon \partial \Delta [n]\xrightarrow{}\Delta [n]$ is
the standard inclusion of simplicial sets, and $\widetilde{f}\boxprod i$
is the pushout product map
\[
(\widetilde{X}\otimes \Delta [n]) \amalg_{\widetilde{X}\otimes \partial
\Delta [n]} (\widetilde{Y}\otimes \partial \Delta [n]) \xrightarrow{}
Y\otimes \Delta [n].
\]
These maps are all cofibrations~\cite[Proposition~5.4.1]{hovey-model}.
Then Proposition~4.2.4 of~\cite{hirschhorn} says that a fibrant object
$X$ is a $\cat{T}$-local fibrant object if and only if
$X\xrightarrow{}*$ has the \rlp the horns on all the maps of $\cat{T}$.

Having obtained the fibrant objects, one defines a map $f$ to be a
\emph{$\cat{T}$-local equivalence} if $\map (f,X)$ is a weak equivalence
of simplicial sets for all $\cat{T}$-local fibrant objects $X$.  These
are the weak equivalences in $L_{\cat{T}}\cat{M}$.  The fibrations are
then the maps that have the \rlp all maps that are both cofibrations and
$\cat{T}$-local equivalences.

If $\cat{M}$ is a simplicial model category, then one can understand the
horns of $f$ by using the simplicial structure.  But unbounded chain
complexes are not simplicial.  It is still easy to understand the horns,
however.  Let $\overline{\Delta}[n]$ be the chain complex of abelian
groups defined by letting $\overline{\Delta [n]}_{k}$ be the free abelian
group on the $\binom{n+1}{k+1}$ $k+1$-element subsets of $\{0,1,\dotsc
,n \}$, for $k\geq 0$.  If
$1_{S}$ denotes the generator corresponding to the set
$S=\{s_{0}<s_{1}<\dotsb < s_{k} \}$, we define 
\[
d(1_{S}) = \sum_{i=0}^{k} (-1)^{i} 1_{S-\{s_{i} \}}. 
\]
This is the obvious chain complex corresponding to the nondegenerate
simplices of $\Delta [n]$.  Then $\overline{\partial \Delta [n]}_{k}$
denotes the subcomplex containing all the $1_{S}$ except the one in
degree $n$ corresponding to $S=\{0,1,\dotsc ,n \}$.  Let
$\overline{i_{n}}$ denote the obvious inclusion $\overline{\partial
\Delta [n]}\xrightarrow{}\overline{\Delta [n]}$.  

The following lemma is a consequence of the naturality of the action of
$\sho \SSet$ on $\sho \cat{M}$, and can be deduced
from~\cite[Chapter~5]{hovey-model}.  

\begin{lemmasec}\label{lem-horns}
Suppose $\cat{M}$ is a $\Ch{\Z}$-model category, and $f\in \cat{M}$ with
a cofibrant approximation $\widetilde{f}$ that is a cofibration.  
Then in the description of Bousfield localization above, one can replace
the horns on $f$ with the maps $\widetilde{f}\boxprod
\overline{i_{n}}$.  
\end{lemmasec}

In general, it is difficult to understand the weak equivalences in
$L_{\cat{T}}\cat{M}$; certainly the maps of $\cat{T}$ become weak
equivalences, but many other maps do as well.  The following theorem is
of some help.

\begin{theoremsec}\label{thm-homotopy-equiv}
Suppose $\cat{M}$ is a left proper model category that is either
cellular or combinatorial, and $\cat{T}$ is a set of maps in $\cat{M}$.
Let $\cat{W}$ be a class of maps in $\cat{M}$ satisfying the two out of
three property, containing the horns on the maps of $\cat{T}$ and every
weak equivalence in $\cat{M}$, and such that maps that are both
cofibrations and in $\cat{W}$ are closed under transfinite compositions
and pushouts.  Then every weak equivalence in $L_{\cat{T}}\cat{M}$ is in
$\cat{W}$.
\end{theoremsec}

\begin{proof}
Suppose $f\mathcolon X\xrightarrow{}Y$ is a weak equivalence in the
Bousfield localization.  Let $L\mathcolon \cat{M}\xrightarrow{}\cat{M}$
denote the functor obtained by applying the small object argument based
on $\cat{T}\cup J$, where $J$ is the set of generating trivial
cofibrations of $\cat{M}$, to the map $X\xrightarrow{}*$.  Then the maps
$X\xrightarrow{}LX$ and $Y\xrightarrow{}LY$ are transfinite compositions
of pushouts of maps of $\cat{T}\cup J$.  This means that they are weak
equivalences in $L_{\cat{T}}\cat{M}$, by Propositions~3.3.10 and
4.2.3 of~\cite{hirschhorn}, and also in $\cat{W}$ by our hypotheses.
Hence $Lf$ is a weak equivalence in $L_{\cat{T}}\cat{M}$ whose domain
and codomain are fibrant (by Proposition~4.2.4 of~\cite{hirschhorn}) in
$L_{\cat{T}}\cat{M}$.  Thus, by Theorem~3.2.13 of~\cite{hirschhorn},
$Lf$ is a weak equivalence in $\cat{M}$ itself, and hence is in
$\cat{W}$.  The two out of three property for $\cat{W}$ now guarantees
that $f$ is in $\cat{W}$.
\end{proof}

In general, Bousfield localization causes one to lose control of the set
of generating trivial cofibrations.  Even if $\cat{M}$ itself has a very
nice set of generating trivial cofibrations, all the theory tells you is
that $L_{\cat{T}}\cat{M}$ has some, possibly gigantic, set of
generating trivial cofibrations.  The following proposition is at least
of some help in dealing with this.

\begin{propositionsec}\label{prop-model-proper}
Suppose $\cat{M}$ is a left proper, cellular or combinatorial model
category, and $\cat{T}$ is a set of maps in $\cat{M}$.  Assume that
$\cat{M}$ has a set of generating trivial cofibrations whose domains are
cofibrant.  Then $L_{\cat{T}}\cat{M}$ has a set of generating trivial
cofibrations whose domains are cofibrant.
\end{propositionsec}

\begin{proof}
Let $J$ be a set of generating trivial cofibrations of $\cat{M}$ whose domains
(and hence codomains) are cofibrant, and let $J'$ be some set of
generating trivial cofibrations of $L_{\cat{T}}\cat{M}$.  For each map
$j\in J'$, choose a cofibration $\hat{j}$ of cofibrant objects that is a
cofibrant approximation to $j$, so that we have a commutative square 
\[
\begin{CD}
\dom \hat{j} @>>> \dom j \\
@V\hat{j}VV @VVjV \\
\codom \hat{j} @>>> \codom j
\end{CD}
\]
where the horizontal maps are weak equivalences in $\cat{M}$.  Let
$\hat{J'}$ denote the set of those $\hat{j}$, and left $K=J\cup
\hat{J'}$.  Then $K$ is a set of trivial cofibrations in
$L_{\cat{T}}\cat{M}$ with cofibrant domains.  We claim that $K$ is a
generating set of trivial cofibrations.  Indeed, suppose $p$ has the
\rlp $K$.  Then $p$ has the \rlp $J$, so $p$ is a fibration in
$\cat{M}$.  Since $p$ also has the \rlp $\hat{J'}$ and $\cat{M}$ is left
proper, Proposition~13.1.16 of~\cite{hirschhorn} implies that $p$ has
the \rlp $J'$, and hence that $p$ is a fibration in
$L_{\cat{T}}\cat{M}$.  
\end{proof}

\section{The homotopy model structure}\label{sec-homotopy}

The object of this section is to construct a model structure on
$\Ch{\Gamma}$, when $\hopfalg$ is an amenable Hopf algebroid, in which
the weak equivalences are the homotopy isomorphisms.
Proposition~\ref{prop-homotopy-projective} tells us that we need to add
more weak equivalences to the projective model structure.  We do this by
using Bousfield localization, described in the previous section.

\subsection{Construction and basic properties}\label{subsec-homotopy-basic}

\begin{definition}\label{defn-homotopy-model}
Suppose $\hopfalg$ is an amenable Hopf algebroid.  Let $\cat{S}$ denote
a set containing one element from each isomorphism class of dualizable
comodules.  Define the \emph{homotopy model structure} on $\Ch{\Gamma}$
to be the Bousfield localization of the projective model structure with
respect to the maps
\[
\eta _{L}\smash S^{n}P\mathcolon S^{n}P \xrightarrow{} LA \smash S^{n}P
\]
for $P\in \cat{S}$ and $n\in \Z$.
\end{definition}

We have already seen that $\Gamma \comod $ is a locally (finitely)
presentable category~\ref{prop-generators}.  It follows easily that
$\Ch{\Gamma }$ is also locally (finitely) presentable, so that
$\Ch{\Gamma }$ is a combinatorial model category.  Thus the
(unpublished) work of Jeff Smith guarantees that the homotopy model
structure exists.  In fact, $\Ch{\Gamma }$ is also cellular, so one can
use Hirschhorn's theory~\cite{hirschhorn} as well.  

Note that the cofibrations do not change under Bousfield localization,
though the fibrations and weak equivalences will change.  This means
that the trivial fibrations also do not change under Bousfield
localization, and therefore that a cofibrant replacement functor in the
projective model structure is also a cofibrant replacement functor in
the homotopy model structure.  Since Bousfield localization preserves
left properness~\cite[Theorem~4.1.1]{hirschhorn}, the homotopy model
structure is left proper.

Our first goal is to prove that the weak equivalences in the homotopy
model structure are the homotopy isomorphisms, explaining the name.

\begin{proposition}\label{prop-homotopy-weak-homotopy}
Let $\hopfalg$ be an amenable Hopf algebroid.  Then every weak
equivalence in the homotopy model structure is a homotopy isomorphism.  
\end{proposition}

\begin{proof}
Proposition~\ref{prop-homotopy-projective} and
Proposition~\ref{prop-homotopy-properties} tell us that the class of
homotopy isomorphisms has all the properties necessary for
Theorem~\ref{thm-homotopy-equiv} to apply.  It remains to check that the
horns on $\eta _{L}\smash S^{n}P$ are homotopy isomorphisms.  Now
$\Ch{\Gamma}$ is a $\Ch{\Z}$ model category; in fact, there is a
symmetric monoidal left Quillen functor
$\Ch{\Z}\xrightarrow{}\Ch{\Gamma}$, induced by the trivial comodule
functor $M \mapsto A\otimes_{\Z}M$.  Lemma~\ref{lem-horns} implies that
the horns on $f$ can be taken to be the maps $f\boxprod (A\otimes_{\Z }
\overline{i_{n}})$.  One can easily check that each map $A\otimes_{\Z }
\overline{i_{n}}$ is a projective cofibration.  The lemma then follows
from Proposition~\ref{prop-homotopy-properties}(e).
\end{proof}

To prove the converse, we need the following proposition.  

\begin{proposition}\label{prop-homotopy-inverts}
Let $\hopfalg$ be a flat Hopf algebroid, and suppose $C$ is cofibrant in
$\Ch{\Gamma}$.  Then $\eta_{L}\smash C\mathcolon
C\xrightarrow{}LA\smash C$ is a weak equivalence in the homotopy model
structure.
\end{proposition}

\begin{proof}
Factor $\eta_{L}\mathcolon S^{0}A\xrightarrow{}LA$ into a cofibration
$i\mathcolon S^{0}A\xrightarrow{}QLA$ followed by a trivial fibration
$q$.  It suffices to show that $i\smash C$ is a trivial cofibration in
the homotopy model structure, because $q\smash C$ is a projective
equivalence by Proposition~\ref{prop-projective-monoid}.  For dualizable
$P$, $\eta _{L}\smash S^{n}P$ is a weak equivalence in the homotopy
model structure by construction, so $i\smash S^{n}P$ is a trivial
cofibration in the homotopy model structure.  Since $C$ is cofibrant,
$0\xrightarrow{}C$ is a retract of a transfinite composition
\[
0=C_{0} \xrightarrow{} C_{1} \xrightarrow{} \dotsb C_{i}\xrightarrow{} \dotsb 
\]
where each map $C_{i}\xrightarrow{}C_{i+1}$ is a pushout of a map
$S^{n-1}P\xrightarrow{}D^{n}P$, where $P\in \cat{S}$ and $n\in \Z$.  It
suffices to show that $i\smash C_{i}$ is a trivial cofibration in the
homotopy model structure for all $i$, which we do by transfinite
induction.  The base case is trivial, since $C_{0}=0$.  For the
successor ordinal step, suppose $i\smash C_{i}$ is a trivial
cofibration.  We have pushout diagrams
\[
\begin{CD}
S^{n-1}P @>>> D^{n}P \\
@VVV @VVV \\
C_{i} @>>> C_{i+1}
\end{CD}
\]
and 
\[
\begin{CD}
QLA\smash S^{n-1}P @>>> QLA\smash D^{n}P \\
@VVV @VVV \\
QLA \smash C_{i} @>>> QLA \smash C_{i+1}
\end{CD}
\]
and $i$ induces a map from one of these to the next.  This map is a weak
equivalence on the upper left corners as mentioned above, a chain
homotopy equivalence on the upper right corners because $D^{n}P$ is
contractible, and a weak equivalence on the lower left corners by the
induction hypothesis.  It follows from the cube
lemma~\cite[Lemma~5.2.6]{hovey-model} that $i\smash C_{i+1}$ is a weak
equivalence as well.  Because the projective structure is monoidal,
$i\smash C_{i+1}$ is also a cofibration.

We are left with the limit ordinal step of the induction.  So suppose
that $i\smash C_{i}$ is a trivial cofibration in the homotopy model
structure for all $i<\alpha $ for some limit ordinal $\alpha$.  Then
Proposition~18.10.1 of~\cite{hirschhorn} implies that $i\smash
C_{\alpha}$ is a weak equivalence as well, and hence a trivial
cofibration.  
\end{proof}

\begin{theorem}\label{thm-homotopy-homotopy}
Suppose $\hopfalg$ is an amenable Hopf algebroid.  Then the weak
equivalences in the homotopy model structure are the homotopy
isomorphisms.  
\end{theorem}
 
\begin{proof}
We have already seen that every weak equivalence in the homotopy model
structure is a homotopy isomorphism.  Conversely, suppose $f\mathcolon
X\xrightarrow{}Y$ is a homotopy isomorphism.  By using a cofibrant
replacement functor $Q$ in the projective model structure, we can
construct the commutative diagram below,
\[
\begin{CD}
QX @>Qf>> QY \\
@Vq_{X}VV @VVq_{Y}V \\
X @>>f> Y
\end{CD}
\]
where $q_{X}$ and $q_{Y}$ are projective equivalences, and $QX$ and $QY$
are cofibrant.  In particular, since projective equivalences are
homotopy isomorphisms by Proposition~\ref{prop-homotopy-projective},
$Qf$ is a homotopy isomorphism.  Since every projective
equivalence is a weak equivalence in the homotopy model structure, it
suffices to show that the homotopy isomorphism $Qf$ is a weak
equivalence.  Consider the commutative square below.
\[
\begin{CD}
QX @>\eta_{L}\smash QX>> LA\smash QX \\
@VQfVV @VVLA\smash QfV \\
QY @>>\eta_{L}\smash QY> LA\smash QY
\end{CD}
\]
Both of the horizontal maps are weak equivalences in the homotopy model
structure by Proposition~\ref{prop-homotopy-inverts}.  Since $Qf$ is a
homotopy isomorphism, $LA\smash Qf$ is a projective equivalence, and
therefore a weak equivalence in the homotopy model structure.  It
follows that $Qf$ is a weak equivalence in the homotopy model
structure as well.  
\end{proof}

Many properties of the homotopy model structure follow immediately from
Theorem~\ref{thm-homotopy-homotopy}.

\begin{theorem}\label{thm-homotopy-model}
Suppose $\hopfalg$ is an amenable Hopf algebroid.  Then the homotopy model
structure is proper, symmetric monoidal, and satisfies the monoid axiom.
Moreover, if $C$ is cofibrant, then $C\smash -$ preserves weak
equivalences.
\end{theorem}

The monoid axiom and the last statement of this theorem are important
because they guarantee the existence of homotopy invariant model
categories of modules and monoids, as explained following the statement
of Proposition~\ref{prop-projective-monoid}.  

\begin{proof}
Bousfield localization preserves left properness, as has already been
observed.  The fact that the homotopy model structure is right proper
follows immediately from Theorem~\ref{thm-homotopy-homotopy} and
Proposition~\ref{prop-homotopy-properties}.  The homotopy model
structure is symmetric monoidal by
Proposition~\ref{prop-homotopy-properties}(e).  

Now for the monoid axiom, which we recall states that any transfinite
composition of pushouts of maps of the form $f\smash X$, where $f$ is a
trivial cofibration and $X$ is arbitrary, is a weak equivalence.  Let us
suppose we know that $f\smash X$ itself is a weak equivalence for all
trivial cofibrations $f$ and all $X$.  Since cofibrations are degreewise
$\Gamma$-split monomorphisms, it follows that $f\smash X$ is also
injective.  Since injective homotopy isomorphisms are closed under
pushouts and filtered colimits by
Propostion~\ref{prop-homotopy-properties}, the monoid axiom will follow.

We are left with showing that $f\smash X$ is a homotopy isomorphism for
all trivial cofibrations $f$ in the homotopy model structure and all
$X$.  It suffices to show this for a set of generating trivial
cofibrations $f$, and Proposition~\ref{prop-model-proper} allows us to
assume those generating trivial cofibrations $f$ have cofibrant domains
and codomains.  Let $q\mathcolon QX\xrightarrow{}X$ be a cofibrant
replacement of $X$, so that $q$ is a projective equvalence and $QX$ is
cofibrant.  We have the commutative square below. 
\[
\begin{CD}
\dom f \smash QX @>\dom f\smash q>> \dom f \smash X \\
@Vf\smash QXVV @VVf\smash XV \\
\codom f \smash QX @>>\codom f\smash q> \codom f \smash X
\end{CD}
\]
By Proposition~\ref{prop-projective-monoid}, both the horizontal maps
are projective equivalences, and hence homotopy isomorphisms.  Since the
homotopy model structure is symmetric monoidal, $f\smash QX$ is a
homotopy isomorphism.  It follows that $f\smash X$ is a homotopy
isomorphism as well, completing the proof of the monoid axiom. 

Now suppose $C$ is cofibrant, and $f$ is a weak equivalence.  Then
$LA\smash f$ is a projective equivalence.  By
Proposition~\ref{prop-projective-monoid}, it follows that $C\smash
LA\smash f$ is still a projective equivalence.  Hence $C\smash f$ is a
homotopy isomorphism, so a weak equivalence.  
\end{proof}

\subsection{Fibrations}\label{subsec-homotopy-fibs}

We now characterize the fibrations in the homotopy model structure.

\begin{proposition}\label{prop-homotopy-fibrations}
Suppose $\hopfalg$ is an amenable Hopf algebroid.  Then a map $p$ is a
fibration in the homotopy model structure if and only if $p$ is a
projective fibration and $\ker p$ is fibrant in the homotopy model
structure.
\end{proposition}

\begin{proof}
Certainly, if $p$ is a homotopy fibration, then $p$ must be a projective
fibration and $\ker p$ must be homotopy fibrant.  Conversely, suppose
$p\mathcolon X\xrightarrow{}Y$ is a projective fibration and $\ker p$ is
homotopy fibrant.  Form the commutative square below by using
factorization,
\begin{equation}\label{eq-square}
\begin{CD}
X @>i_{X}>> X' \\
@VpVV @VVqV \\
Y @>>i_{Y}> Y'
\end{CD}
\end{equation}
where $i_{X}$ and $i_{Y}$ are trivial cofibrations in the homotopy model
structure, $X'$ and $Y'$ are homotopy fibrant, and $q$ is a homotopy
fibration.  We claim that this square is a homotopy pullback square in
the projective model structure.  Proposition~3.4.7 of~\cite{hirschhorn}
then implies that $p$ is a fibration in the homotopy model structure.

To see that the square~\ref{eq-square} is a homotopy pullback square,
let $P\xrightarrow{q'}Y$ be the pullback of $q$ through $i_{Y}$.  Then
$q'$ is a projective fibration, and there is an induced factorization
\[
X\xrightarrow{s} P \xrightarrow{t} X'
\]
of $i_{X}$.  Since $t$ is the pullback of the homotopy isomorphism
$i_{Y}$ through the surjection $q$,
Proposition~\ref{prop-homotopy-properties} implies that $t$ is a
homotopy isomorphism.  Hence $s$ is a homotopy isomorphism as well.
Consider the commutative diagram below
\[
\begin{CD}
0 @>>> \ker p @>>> X @>p>> Y @>>> 0 \\
@. @VrVV @VsVV @| \\
0 @>>> \ker q @>>> P @>q'>> Y @>>> 0
\end{CD}
\]
whose rows are short exact (since projective fibrations are surjective).
The long exact sequence in homotopy implies that $r$ is a homotopy
isomorphism.  But $\ker p$ is homotopy fibrant by assumption, and $\ker
q$ is homotopy fibrant since $q$ is a homotopy fibration.
Theorem~3.2.13 of~\cite{hirschhorn} implies that $r$ is a projective
equivalence.  Applying $\Gamma \comod (Q,-)$ for $Q\in \cat{S}$ and
considering the long exact homology sequence shows that $s$ is a
projective equivalence as well.  This means that the
square~\ref{eq-square} is a homotopy pullback square, completing the
proof.
\end{proof}

The characterization of fibrations we have just given would be more
helpful if we knew what the fibrant objects in the homotopy model
structure are.  

\begin{theorem}\label{thm-homotopy-fibrant}
Suppose $\hopfalg$ is an amenable Hopf algebroid.  Then the following
are equivalent.
\begin{enumerate}
\item [(a)] $\eta_{L}\smash X\mathcolon X\xrightarrow{}LA\smash X$ is a
projective equivalence. 
\item [(b)] $X$ is projectively equivalent to some complex of relative
injectives.  
\item [(c)] $X$ is fibrant in the homotopy model structure.  
\end{enumerate}
\end{theorem}

\begin{proof}
It is clear that (a) implies (b).  To see that~(b) implies~(c), our
first goal is to show that if $X$ is projectively equivalent to some
complex of relative injectives, then $\eta _{L}\smash X$ is a projective
equivalence (incidentally proving (b) implies (a)).  It obviously
suffices to show this for actual complexes of relative injectives $X$.
Any such complex can be written as the colimit of the bounded above
complexes $X^{n}$, where $X^{n}_{i}=0$ for $i>n$ and $X^{n}_{i}=X_{i}$
for $i\leq n$.  Since colimits of projective equivalences are projective
equivalences, we can assume that $X$ is a bounded above complex of
relative injectives.  In this case, we will show that $\eta _{L}\smash
X$ is in fact a chain homotopy equivalence.  Indeed, since $\eta_{L}$ is
a degreewise $A$-split monomorphism, $\eta_{L}\smash X$ is a degreewise
$A$-split monomorphism between complexes of relative injectives.  It
follows that $\eta_{L}\smash X$ is a degreewise split monomorphism of
relative injectives.  Let $Y$ denote the cokernel of $\eta_{L}\smash X$.
Then $Y$ is also a bounded above complex of relative injectives, and $Y$
is contractible as a complex of $A$-modules since $\eta_{L}$ is a chain
homotopy equivalence of complexes of $A$-modules.
Lemma~\ref{lem-contractible-to-injective} then implies that $Y$ is
contractible as a complex of comodules.  Given this, an elementary
argument then shows that $\eta_{L}\smash X$ is a chain homotopy
equivalence.

Now suppose that $X$ is projectively equivalent to a complex of relative
injectives, $C$ is cofibrant and $\pi_{*}C=0$.  We claim
that every chain map $f\mathcolon C\xrightarrow{}X$ is chain homotopic
to $0$.  Indeed, the composite
\[
C \xrightarrow{f} X \xrightarrow{\eta_{L}\smash X} LA\smash X
\]
factors through $LA\smash C$, which is projectively trivial.  Hence the
map $(\eta_{L}\smash X)\circ f$ is $0$ in the homotopy category of the
projective model structure.  Since $\eta _{L}\smash X$ is a projective
equivalence by the previous paragraph, we conclude that $f$ is $0$ in
the homotopy category of the projective model structure.  Since $C$ is
cofibrant and everything is fibrant in the projective model structure,
it follows that $f$ is chain homotopic to $0$.

We can now complete the proof that (b) implies (c) by showing that
$X\xrightarrow{}0$ has the \rlp every cofibration $i\mathcolon
A\xrightarrow{}B$ that is a homotopy isomorphism.  Indeed, suppose
$f\mathcolon A\xrightarrow{}X$ and let $C$ be the cokernel of $i$.
Since $i$ is a cofibration, it is a split monomorphism in each degree,
and so $B_{n}\cong A_{n}\oplus C_{n}$.  The differential on $B$ must
then be given by $d(a,c)=(da+\tau c,dc)$, where $\tau d=-d\tau$.  Thus
$\tau$ is a chain map from the desuspension $\Sigma^{-1}C$ of $C$ to
$A$.  The composition
\[
\Sigma^{-1} C \xrightarrow{\tau} A \xrightarrow{f} X
\]
must be chain homotopic to $0$, since $\Sigma^{-1}C$ is cofibrant and
$\pi_{*}\Sigma^{-1}C=0$.  Hence there are maps $D_{n}\mathcolon
C_{n}\xrightarrow{}A_{n}$ such that $-D_{n-1}d+dD_{n}=f\tau $.  Define
\[
g(a,c)=fa+D_{n}c.
\]
Then $g\mathcolon B\xrightarrow{}C$ is a chain map extending $f$, so
$X\xrightarrow{}0$ has the \rlp $i$ as required.

We now show that (c) implies (a).  So suppose $X$ is fibrant in the
homotopy model structure.  The map $\eta_{L}\smash X\mathcolon
X\xrightarrow{}LA\smash X$ is a homotopy isomorphism by
Corollary~\ref{cor-homotopy-none}.  Since $X$ is fibrant, and $LA\smash
X$ is also fibrant since~(b) implies~(c), it follows from Theorem~3.2.13
of~\cite{hirschhorn} that $\eta_{L}\smash X$ is a projective
equivalence. 
\end{proof}

\begin{corollary}\label{cor-homotopy-fibrant}
Suppose $\hopfalg$ is an amenable Hopf algebroid, and give $\Ch{\Gamma}$
the homotopy model structure.  For any $X\in \Ch{\Gamma}$, the map 
\[
\eta_{L}\smash X\mathcolon X\xrightarrow{}LA\smash X
\]
is a natural weak equivalence whose target is fibrant.
\end{corollary}

This follows immediately from Corollary~\ref{cor-homotopy-none} and
Theorem~\ref{thm-homotopy-fibrant}.  Note that $\eta_{L}\smash X$ is not
normally a cofibration, however.

We also note the following corollary.

\begin{corollary}\label{cor-amenable-fin-gen}
Suppose $\hopfalg $ is an amenable Hopf algebroid, and given $\Ch{\Gamma
}$ the homotopy model structure.  Then weak equivalences and fibrations
are closed under filtered colimits.  
\end{corollary}

This corollary is saying that the homotopy model structure behaves as if
it were finitely generated.  We do not know if it is in fact finitely
generated for a general amenable Hopf algebroid, though for an Adams
Hopf algebroid it is.

\begin{proof}
We have seen in Proposition~\ref{prop-homotopy-properties} that homotopy
isomorphisms are closed under filtered colimits.  Since homotopy
fibrations are just projective fibrations with homotopy fibrant kernel,
and projective fibrations are closed under filtered colimits, it
suffices to show that homotopy fibrant objects are closed under filtered
colimits.  This follows from the characterization of homotopy fibrant
objects in part~(a) of Theorem~\ref{thm-homotopy-fibrant}, since
projective equivalences are closed under filtered colimits.  
\end{proof}

\subsection{Naturality}\label{subsec-homotopy-natural}

Like the projective model structure, the homotopy model structure is
natural.

\begin{proposition}\label{prop-homotopy-natural}
Suppose $\Phi \mathcolon \hopfalg \xrightarrow{}\otherhopfalg$ is a map
of amenable Hopf algebroids.  Then $\Phi$ induces a left Quillen functor
$\Phi_{*}\mathcolon \Ch{\Gamma}\xrightarrow{}\Ch{\Sigma}$ of the
homotopy model structures.  
\end{proposition}

\begin{proof}
By Proposition~\ref{prop-projective-natural}, $\Phi^{*}$ is a right
Quillen functor of the projective model structures.  Thus $\Phi ^{*}$
preserves projective fibrations and projective equivalences (since
everything is fibrant in the projective model structure), and so will
also preserves trivial fibrations in the homotopy model structure, since
these coincide with projective trivial fibrations.  Suppose $p$ is a
homotopy fibration.  Then $p$ is a projective fibration such that $\ker
p$ is projectively equivalent to a complex of relative injectives $K$,
by Proposition~\ref{prop-homotopy-fibrations} and
Theorem~\ref{thm-homotopy-fibrant}.  Hence $\Phi ^{*}p$ is also a
projective fibration, and $\ker \Phi ^{*}p$ is projectively equivalent
to $\Phi ^{*}K$.  We will show that $\Phi ^{*}K$ is a complex of
relative injectives.  Hence $\Phi ^{*}p$ is a homotopy fibration by
Theorem~\ref{thm-homotopy-fibrant} and
Proposition~\ref{prop-homotopy-fibrations}, and so $\Phi ^{*}$ is a
right Quillen functor as required.

We are now reduced to showing that $\Phi ^{*}$ preserves relative
injectives.  Suppose $I$ is a relatively injective $\Sigma $-comodule,
and $E$ is an $A$-split short exact sequence of $\Gamma $-comodules.
Then 
\[
\Gamma \comod (E, \Phi ^{*}I) \cong \Sigma \comod (\Phi _{*}E,I).
\]
Since $\Phi _{*}E=B\otimes _{A}E$, $\Phi _{*}E$ is a $B$-split short
exact sequence, so $\Sigma \comod (\Phi _{*}E,I)$ is exact.  
\end{proof}

The homotopy model structure is also invariant under weak equivalences,
but this is considerably harder to prove.  We begin with a definition. 

\begin{definition}\label{defn-pseudo-injective}
Suppose $\hopfalg$ is an amenable Hopf algebroid.  Define a
$\Gamma$-comodule $I$ to be \emph{pseudo-injective} if
$\Ext^{n}_{\Gamma}(P,I)=0$ for all dualizable comodules $P$.  
\end{definition}

Every relative injective is pseudo-injective, by
Lemma~\ref{lem-proj-inj}.  The reason for introducing pseudo-injectives
is the following lemma, which would be false for relative injectives.  

\begin{lemma}\label{lem-homotopy-weak-pseudo}
Suppose $\Phi \mathcolon \hopfalg \xrightarrow{}\otherhopfalg$ is a weak
equivalence of flat Hopf algebroids.  If $I$ is a pseudo-injective
$\Gamma$-comodule, then $\Phi_{*}I$ is a pseudo-injective
$\Sigma$-comodule. 
\end{lemma}

\begin{proof}
Suppose $P$ is a dualizable $\Sigma$-comodule.  Because $\Phi^{*}$ is an
equivalence of categories whose right adjoint is naturally isomorphic to
$\Phi_{*}$, we have
\[
\Ext^{n}_{\Sigma}(P, \Phi_{*}I) \cong \Ext^{n}_{\Gamma}(\Phi^{*}P, I).
\]
As explained in the proof of Theorem~\ref{thm-proj-weak}, $\Phi ^{*}P$
is a dualizable $\Gamma $-comodule.  The lemma follows.
\end{proof}

\begin{lemma}\label{lem-homotopy-pseudo-good}
Suppose $\hopfalg$ is a flat Hopf algebroid.  Any bounded above complex
of pseudo-injectives with no homology is projectively trivial.  
\end{lemma}

\begin{proof}
Suppose $X$ is a bounded above complex of pseudo-injectives with no
homology.  We claim that the cycle comodule $Z_{n}X$ is pseudo-injective
for all $n$.  This is obvious for large $n$, since $X$ is bounded
above.  We have a short exact sequence 
\[
0 \xrightarrow{} Z_{n}X \xrightarrow{} X_{n} \xrightarrow{} Z_{n-1}X
\xrightarrow{} 0
\]
since $X$ has no homology.  The long exact sequence in
$\Ext_{\Gamma}^{*}(P,-)$ shows that $Z_{n-1}X$ is pseudo-injective.  By
induction, $Z_{n}X$ is pseudo-injective for all $n$.  One can then
easily check that $\Gamma \comod (P,X)$ is exact for all dualizable
comodules $P$.
\end{proof}

\begin{corollary}\label{cor-homotopy-pseudo-homology}
Suppose $\hopfalg$ is a flat Hopf algebroid, and $f\mathcolon
X\xrightarrow{}Y$ is a homology isomorphism of complexes of bounded
above pseudo-injectives.  Then $f$ is a projective equivalence. 
\end{corollary}

\begin{proof}
Let $C$ denote the mapping cylinder of $f$, and let $Z=C/X$ denote the
mapping cone of $f$.  Then $f=pi$, where $i\mathcolon X\xrightarrow{}C$
is a degreewise split monomorphism and $p$ is a chain homotopy
equivalence.  It therefore suffices to show that the homology
isomorphism $i$ is a projective
equivalence.  Since $i$ is degreewise split, it suffices to show that
$Z$ is projectively trivial.  But $Z_{n}=Y_{n}\oplus 
X_{n-1}$, so $Z$ is a bounded above complex of pseudo-injectives with
no homology.  Lemma~\ref{lem-homotopy-pseudo-good} implies that $Z$ is
projectively trivial.  
\end{proof}

\begin{theorem}\label{thm-homotopy-natural-weak}
Suppose $\Phi \mathcolon \hopfalg \xrightarrow{}\otherhopfalg$ is a weak
equivalence of Adams Hopf algebroids.  Then $\Phi_{*}\mathcolon
\Ch{\Gamma}\xrightarrow{}\Ch{\Sigma}$ is a Quillen equivalence of the
homotopy model structures.   In fact, both $\Phi_{*}$ and $\Phi^{*}$
preserve and reflect homotopy isomorphisms.  
\end{theorem}

\begin{proof}
We have seen in Theorem~\ref{thm-proj-weak} that $\Phi_{*}$ is a
Quillen equivalence of the projective model structures, and that $\Phi
_{*}$ and $\Phi ^{*}$ preserve and reflect projective equivalences.  We
first show that $\Phi _{*}$ preserves homotopy isomorphisms, and hence
that $\Phi ^{*}$ reflects homotopy isomorphisms.  Indeed, it follows
from Proposition~\ref{prop-homotopy-natural} that $\Phi _{*}$ preserves
homotopy trivial cofibrations.  Thus it suffices to show that $\Phi
_{*}p$ is a homotopy isomorphism when $p$ is a homotopy trivial
fibration.  But then $p$ is a projective equivalence, so $\Phi _{*}p$ is
also a projective equivalence.

It is more difficult to show that $\Phi _{*}$ reflects homotopy
isomorphisms.  To see this, we first construct a factorization 
\[
B=\Phi_{*}A \xrightarrow{\Phi_{*}\eta_{L}} \Phi_{*}LA
\xrightarrow{\alpha} LB
\]
of $\eta_{L}\mathcolon B\xrightarrow{}LB$.  We have
\[
(\Phi_{*}LA)_{n}\cong \Phi_{*}\Gamma \smash (\Phi_{*}\overline{\Gamma
})^{\smash n} \text{ and } (LB)_{n}=\Sigma \smash
\overline{\Sigma}^{\smash n}.
\]
There is a natural map of comodules $\Phi_{*}\Gamma =B\otimes
\Gamma \xrightarrow{}\Sigma$ that takes $b\otimes x$ to
$\eta_{L}(b)\Phi_{1}(x)$.  This map induces a map
$\Phi_{*}\overline{\Gamma}\xrightarrow{}\overline{\Sigma}$, which in
turn induces the desired map $\alpha \mathcolon
\Phi_{*}LA\xrightarrow{}LB$ of complexes.

Now $\Phi_{*}LA$ and $LB$ are both complexes of pseudo-injectives,
by Lemma~\ref{lem-homotopy-weak-pseudo}.  The map $\Phi_{*}\eta_{L}$ is
a homology isomorphism, since $\Phi_{*}$, like any equivalence of
abelian categories, is exact.  The map $\eta_{L}\mathcolon
B\xrightarrow{}LB$ is a homology isomorphism by construction.  Thus,
$\alpha \mathcolon \Phi_{*}LA\xrightarrow{}LB$ is a homology isomorphism
of bounded above complexes of pseudo-injectives, and so a projective
equivalence, by Corollary~\ref{cor-homotopy-pseudo-homology}.

We can now show that $\Phi _{*}$ reflects homotopy isomorphisms between
cofibrant objects.  Indeed, suppose that $f\mathcolon X\xrightarrow{}Y$
is a map of cofibrant objects such that $\Phi_{*}f$ is a homotopy
isomorphism.  Then $LB\smash \Phi_{*}f$ is a projective equivalence.
But we have the commutative square below.
\[
\begin{CD}
\Phi_{*}LA \smash \Phi_{*}X @>\Phi_{*}LA\smash \Phi_{*}f>> \Phi_{*}LA
\smash \Phi_{*}Y \\ 
@V\alpha \smash \Phi_{*}XVV @VV\alpha \smash \Phi_{*}YV \\
LB \smash \Phi_{*}X @>>LB\smash \Phi_{*}f> LB\smash \Phi_{*}Y
\end{CD}
\]
Proposition~\ref{prop-projective-monoid} implies that the vertical maps
in this square are projective equivalences.  Hence $\Phi_{*}LA\smash
\Phi_{*}f\cong \Phi_{*}(LA\smash f)$ is a projective equivalence.  But
$\Phi _{*}$ reflects projective equivalences, so $LA\smash f$ is a
projective equivalence.  Hence $f$ is a homotopy isomorphism as
required.  

We now claim that $\Phi_{*}$ reflects all homotopy isomorphisms, from
which it follows easily that $\Phi ^{*}$ preserves all homotopy
isomorphisms.  So suppose $f\mathcolon X\xrightarrow{}Y$ is a map such
that $\Phi_{*}f$ is a homotopy isomorphism.  We have the commutative
square below, in which the vertical maps are projective equivalences and
$QX$ and $QY$ are cofibrant. 
\[
\begin{CD}
QX @>Qf>> QY \\
@Vq_{X}VV  @VVq_{Y}V \\
X @>>f> Y
\end{CD}
\]
Since $\Phi_{*}$ preserves projective equivalences, we conclude that
$\Phi_{*}Qf$ is a homotopy isomorphism.  Since $\Phi_{*}$ reflects
homotopy isomorphisms between cofibrant objects, $Qf$ is a homotopy
isomorphism, and hence $f$ is a homotopy isomorphism as required.  

It now follows easily that $\Phi _{*}$ is a Quillen equivalence, as in
the proof of Theorem~\ref{thm-proj-weak}.  
\end{proof}

\subsection{Comparison with injective model
structure}\label{subsec-homotopy-comp} 

When $A=k$ is a field, there is a model structure on $\Ch{\Gamma}$ in
which $\sho \Ch{\Gamma}(A,A)_{*}\cong \Ext^{-*}_{\Gamma}(A,A)$ developed
in~\cite[Section~2.5]{hovey-model}.  In this model structure, which we
call the \emph{injective model structure}, the cofibrations are just the
monomorphisms, and the fibrations are the degreewise surjections with
degreewise injective kernels.  (Remember that relatively injective and
injective coincide in case $A$ is a field).  The weak equivalences are
the homotopy isomorphisms, where homotopy is defined as in
Definition~\ref{defn-homotopy} but using only simple comodules (that is,
those comodules with no nontrivial subcomodule) as the source.  For
years, the author searched for a generalization of this model structure
to Hopf algebroids without success.  The injective model structure is
\textbf{NOT} a special case of the homotopy model structure.  Indeed,
cofibrations in the homotopy model structure are degreewise split over
$\Gamma$, whereas cofibrations in the injective model structure are
split over $A$, but not necessarily $\Gamma$.  However, we have the
following theorem.

\begin{theorem}\label{thm-homotopy-comparison}
Suppose $\Gamma $ is a Hopf algebra over a field $k=A$.  Then the
identity functor defines a Quillen equivalence from the homotopy model
structure to the injective model structure.   
\end{theorem}

\begin{proof}
We claim that the two model structures have the same weak equivalences.
If we can prove this, then the identity functor will be a left Quillen
functor from the homotopy model structure to the injective model
structure, since every projective cofibration is a monomorphism.  It
must be a Quillen equivalence since the weak equivalences are the same.  

Note that the dualizable comodules coincide with the finite-dimensional
comodules.  Every simple comodule is finite-dimensional by Lemma~9.5.5
of~\cite{hovey-axiomatic}.  Thus every homotopy isomorphism is a weak
equivalence in the injective model structure.  To prove the converse, it
suffices to prove that if $\Gamma \comod (P, LA \smash f)$ is a homology
isomorphism for all simple comodules $P$, then it is a homology
isomorphism for all finite-dimensional comodules $P$.  We do this by
induction on the dimension of $P$.  Since every one-dimensional comodule
is simple, the base case is easy.  Now suppose we know $\Gamma \comod
(P,LA\smash f)$ is an isomorphism for all $P$ of dimension $<n$, and $P$
has dimension $n$.  If $P$ is simple, there is nothing to prove.  If $P$
is not simple, there is a short exact sequence of comodules
\[
0 \xrightarrow{} Q \xrightarrow{} P \xrightarrow{} P/Q \xrightarrow{} 0
\]
with $\dim Q<n$.  This sequence is necessarily split over $k$, since $k$
is a field.  Therefore, the sequence 
\[
0 \xrightarrow{}\Gamma \comod (P/Q,LA\smash X) \xrightarrow{} \Gamma
\comod (P,LA\smash X) \xrightarrow{} \Gamma \comod (Q,LA\smash X)
\xrightarrow{} 0
\]
is still short exact, as is the corresponding sequence with $Y$
replacing $X$.  The map $f$ induces a map between the corresponding long
exact sequences in homology.   By the induction hypothesis, $\Gamma
\comod (Q,LA\smash f)$ and $\Gamma \comod (P/Q,LA\smash f)$ are
isommorphism.  The five lemma implies that $\Gamma \comod (P,LA\smash
f)$ is an isomorphism as well, completing the proof of the induction
step.  
\end{proof}

\section{The stable category}\label{sec-stable}

We define the homotopy category of the homotopy model structure on
$\Ch{\Gamma }$ to be the \emph{stable homotopy category of $\hopfalg $},
and we denote it by $\stable{\Gamma }$, following
Palmieri~\cite{palmieri-book} in the case of the Steenrod algebra.  The
category $\stable{\Gamma }$ is what we should mean by the \emph{derived
category} $\cat{D}\hopfalg $ of the Hopf algebroid $(A, \Gamma )$.
This is consistent with the usual notation, since
$\cat{D}(A,A)=\cat{D}(A)$, the usual unbounded derived category of $A$.
However, we must remember that to form the derived category, we invert
the \textbf{homotopy isomorphisms}, not the homology isomorphisms.  It
is just that in the case of a discrete Hopf algebroid $(A,A)$, the
homotopy isomorphisms coincide with the homology isomorphisms.  

We conclude the paper by establishing some basic properties of
$\stable{\Gamma }$.  We show that it is a unital algebraic stable
homotopy category~\cite{hovey-axiomatic}.  This means that it shares
most of the formal properties of the derived category of a commutative
ring, or the ordinary stable homotopy category, except that it has
several generators rather than just one.  In certain cases of interest
in algebraic topology, such as $\Gamma =BP_{*}BP$, we show that
$\stable{\Gamma }$ is monogenic, so that (bigraded) suspensions of
$BP_{*}$ weakly generate the category.  We also show that
\[
\stable{\Gamma }(S^{0}M, S^{k}N)\cong \Ext _{\Gamma }^{k}(M,N)
\]
for certain $\Gamma $-comodules $M$ and $N$. 

We begin with the following lemma.

\begin{lemmasec}\label{lem-strong-dualizable}
Suppose $\hopfalg$ is an amenable Hopf algebroid, and $P$ is a
dualizable $\Gamma$-comodule.  Then $S^{n}P$ is dualizable in the
homotopy category of the projective model structure on $\Ch{\Gamma}$ for
all $n$.
\end{lemmasec}

\begin{proof}
Recall that the symmetric monoidal product in the homotopy category is
the derived smash product $X\smash ^{L}Y=QX\smash QY$, where $Q$ denotes
a cofibrant replacement functor.  Similarly, the closed structure is
$RF(X,Y)=F(QX,RY)$, where $R$ is a fibrant replacement functor.  To show
that $X$ is dualizable, we must show that the unit
\[
S^{0}A\xrightarrow{}RF(X,X)
\]
factors through the composition map
\[
RF(X,S^{0}A)\smash ^{L}X\xrightarrow{}RF(X,X).
\]
In the projective model structure, everything is fibrant, so we may as
well take $R$ to be the identity functor.  Furthermore, $S^{n}P$ is
cofibrant, so we conclude that 
\[
RF(S^{n}P, S^{n}P)\cong S^{0}F(P, P), RF(S^{n}P,S^{0}A) \cong
S^{-n}F(P,A), 
\]
and 
\[
RF(S^{n}P. S^{0}A) \smash ^{L}S^{n}P \cong S^{0}(F(P,A) \smash P).
\]
It is now clear that $S^{n}P$ is dualizable, since $P$ is so.  
\end{proof}

\begin{theoremsec}\label{thm-stable-homotopy-cat}
Suppose $\hopfalg$ is an amenable Hopf algebroid.  Then the homotopy
category of the projective model structure and $\stable{\Gamma }$ are unital
algebraic stable homotopy categories.  A set of small, dualizable, weak
generators is given by the set of all $S^{n}P$ for $P$ a dualizable
comodule and $n\in \Z$.
\end{theoremsec}

\begin{proof}
It is easy to check that the ordinary suspension, defined by $(\Sigma
X)_{n}=X_{n-1}$ with $d_{\Sigma X}=-d_{X}$, is a Quillen equivalence of
both the projective and homotopy model structures.  One can also check
that it induces the model category theoretic suspension on the homotopy
categories.  This means that both the projective and homotopy model
structures are stable in the sense of~\cite[Section~7.1]{hovey-model},
and therefore that the homotopy category of the projective model
structure and $\stable{\Gamma }$ are triangulated.

Since the projective and homotopy model structures are symmetric
monoidal, their homotopy categories are also symmetric monoidal in a way
that is compatible with the triangulation (see Chapters~4 and~6
of~\cite{hovey-model}).  In fact, they satisfy much stronger
compatibility relations than those demanded in~\cite{hovey-axiomatic};
see~\cite{may-triangulated}.

The projective model structure is finitely generated, so the results of
Sections~7.3 and~7.4 of~\cite{hovey-model} guarantee that the cofibers
of the generating cofibrations, namely the $S^{n}P$, form a set of small
weak generators for the homotopy category.  The homotopy model structure
may not be finitely generated, but fibrations and weak equivalences are
closed under filtered colimits by Corollary~\ref{cor-amenable-fin-gen},
and this is all that is needed for the arguments of Section~7.4
of~\cite{hovey-model} to apply.  Thus the $S^{n}P$ also form a set of
small weak generators for $\stable{\Gamma }$.
Lemma~\ref{lem-strong-dualizable} guarantees that they are dualizable in
the homotopy category of the projective model structure; since the
functor from this category to $\stable{\Gamma }$ is symmetric monoidal,
they are also dualizable in $\stable{\Gamma }$.
\end{proof}

We now investigate when the homotopy category of the homotopy model
structure is monogenic.  We first recall a definition from abelian
categories.

\begin{definitionsec}\label{defn-wide}
A \emph{thick subcategory} of an abelian category $\cat{C}$ is 
a full subcategory $\cat{T}$ that is closed under retracts and has the
two-out-of-three property.  This means that if 
\[
0 \xrightarrow{} M' \xrightarrow{} M \xrightarrow{} M'' \xrightarrow{} 0
\]
is a short exact sequence, and two out of $M', M, M''$ are in $\cat{T}$,
so is the third.  If $\cat{C}$ is graded, we also insist that thick
subcategories be closed under aribtrary shifts.
\end{definitionsec}

The reader used to $BP_{*}BP$-comodules will be familiar with the
following definition.  

\begin{definitionsec}\label{defn-Landweber-filtration}
Suppose $\hopfalg $ is an amenable Hopf algebroid, and $M$ is a $\Gamma
$-comodule.  We say that \emph{$M$ has a Landweber filtration} if there
is a finite filtration 
\[
0=M_{0} \subseteq M_{1}\subseteq \dots \subseteq M_{t}=M
\]
of $M$ by subcomodules such that each quotient $M_{j}/M_{j-1}\cong
A/I_{j}$ for some ideal $I_{j}$ of $A$ that is generated by an
invariant finite regular sequence.  
\end{definitionsec}

We recall that a sequence $x_{1},\dots ,x_{n}$ is an invariant regular
sequence if $x_{i}$ is a primitive nonzero divisor in $A/(x_{1},\dots
,x_{i-1})$ for all $i$.

In case $\hopfalg $ is graded, we allow the filtration quotients
$M_{i}/M_{i-1}$ to be isomorphic to some shift of $A/I_{j}$ rather than
$A/I_{j}$ itself.  

From a structural point of view, whether or not $M$ has a Landweber
filtration is not important.  What matters is whether $M$ is in the thick
subcategory generated by $A$.  

\begin{lemmasec}\label{lem-wide}
Suppose $\hopfalg $ is an amenable Hopf algebroid, and $M$ is a $\Gamma
$-comodule with a Landweber filtration.  Then $M$ is in the thick
subcategory generated by $A$.  
\end{lemmasec}

\begin{proof}
Since thick subcategories are closed under extensions, it suffices to
check that $A/I$ is in the thick subcategory generated by $A$, where $I$
is generated by a finite invariant regular sequence $x_{1},\dots
,x_{n}$.  This follows from the short exact sequences of comodules
\[
0\xrightarrow{}A/(x_{1},\dots ,x_{i-1}) \xrightarrow{x_{i}}
A/(x_{1},\dots ,x_{i-1}) \xrightarrow{} A/(x_{1},\dots ,x_{i})
\xrightarrow{} 0
\]
and induction.  
\end{proof}

\begin{theoremsec}\label{thm-stable-monogenic}
Suppose $\hopfalg$ is an amenable Hopf algebroid, and every dualizable
comodule $P$ is in the thick subcategory generated by $A$.  Then
$\stable{\Gamma }$ is monogenic, in the sense that $\{S^{n}A \}$ form a
set of small weak generators.
\end{theoremsec}

In the graded case, we would instead get that $\{S^{n,m}A \}$ would form
a set of weak generators.  

\begin{corollarysec}\label{cor-stable-monogenic}
Let $E$ be a ring spectrum that is Landweber exact over $MU$ or $BPJ$
for some finite invariant regular sequence $J$, and suppose that
$E_{*}E$ is commutative.  Then $\stable{E_{*}E}$ is a bigraded monogenic
stable homotopy category.
\end{corollarysec}

\begin{proof}
It is shown in~\cite{hovey-strickland-comodules} that every finitely
presented $E_{*}E$-comodule is a retract of a comodule with a Landweber
filtration, and hence in the thick subcategory generated by $E_{*}$.  
\end{proof}

To prove Theorem~\ref{thm-stable-monogenic}, we first need a lemma.

\begin{lemmasec}\label{lem-stable-cofiber}
Suppose $\hopfalg $ is an amenable Hopf algebroid, and 
\[
0 \xrightarrow{} M' \xrightarrow{f} M \xrightarrow{} M'' \xrightarrow{} 0
\]
is a short exact sequence of comodules.  Then 
\[
S^{0}M' \xrightarrow{} S^{0}M \xrightarrow{} S^{0}M''  
\]
is a cofiber sequence in $\stable{\Gamma }$.  
\end{lemmasec}

\begin{proof}
Factor $S^{0}f$ into a projective cofibration $i\mathcolon
S^{0}M'\xrightarrow{}X$ followed by a projective trivial fibration $p$.
Then we have the commutative diagram below, whose rows are exact. 
\[
\begin{CD}
0 @>>> S^{0}M' @>i>> S^{0}M @>>> C @>>> 0 \\
@. @| @VpVV @VVqV \\
0 @>>> S^{0}M' @>>S^{0}f> S^{0}M @>>> S^{0}M'' @>>> 0 \\
\end{CD}
\]
The long exact sequence in homotopy
(Lemma~\ref{lem-homotopy-properties}) and the five lemma imply that $q$
is a homotopy isomorphism.  Therefore the bottom row is isomorphic in
$\stable{\Gamma }$ to the top row,
which is a cofiber sequence.
\end{proof}

\begin{proof}[Proof of Theorem~\ref{thm-stable-monogenic}]
Suppose that $\pi^{A}_{*}(X)=0$.  Let us denote maps in the homotopy
category of the homotopy model structure by $[Y,Z]_{*}$, so that
$[S^{0}A,X]_{*}=0$.   Let $\cat{T}$ denote the full subcategory of all
comodules $M$ such that $[S^{0}M, X]_{*}=0$.  We claim that $\cat{T}$ is a
thick subcategory, and therefore contains the dualizable comodules.
Theorem~\ref{thm-stable-homotopy-cat} then completes the proof. 

It is clear that $\cat{T}$ is closed under retracts.  To show that
$\cat{T}$ is thick, suppose we have a short exact sequence 
\[
0 \xrightarrow{} M' \xrightarrow{} M \xrightarrow{} M'' \xrightarrow{} 0
\] 
such that two out of $M', M, M''$ are in $\cat{T}$.
Lemma~\ref{lem-stable-cofiber} then implies that 
\[
S^{0}M' \xrightarrow{} S^{0}M \xrightarrow{} S^{0}M''
\]
is a cofiber sequence in $\stable{\Gamma }$.  The long exact sequence
obtained by applying $[, X]_{*}$ then shows that the other one is also
in $\cat{T}$.  
\end{proof}

Finally, we study $\stable{\Gamma }(S^{0}M,S^{0}N)$ for comodules $M$
and $N$.  

\begin{definitionsec}\label{defn-localizing}
A full subcategory of an abelian category $\cat{C}$ is called
\emph{localizing} if it is a thick subcategory closed under coproducts.  
\end{definitionsec}

\begin{propositionsec}\label{prop-ext}
Let $\hopfalg $ be an amenable Hopf algebroid, and $M$ and $N$ be
$\Gamma $-comodules.  Then there is a natural map 
\[
\Ext ^{k}_{\Gamma }(M,N) \xrightarrow{\alpha _{MN}} \stable{\Gamma
}(S^{0}M, S^{k}N).
\]
This map is an isomorphism if $M$ is in the localizing subcategory
generated by the dualizable comodules.  
\end{propositionsec}

Note that the $\Ext $ groups that appear in this proposition are
$\Ext $ groups in the category of $\Gamma $-comodules, not relative
$\Ext $ groups.  

\begin{proof}
A class in $\Ext ^{k}_{\Gamma }(M,N)$ is represented by an exact
sequence of comodules
\[
0 \xrightarrow{} N=E_{0} \xrightarrow{f_{0}} E_{1} \xrightarrow{f_{1}}
E_{2}\xrightarrow{f_{2}}\dots \xrightarrow{f_{k-1}} E_{k}
\xrightarrow{f_{k}} E_{k+1}=M \xrightarrow{} 0.
\]
We can split this into the short exact sequences 
\[
0 \xrightarrow{} \ker f_{i} \xrightarrow{} E_{i} \xrightarrow{} \coker
f_{i} \xrightarrow{} 0. 
\]
Each such short exact sequence is gives rise to a cofiber sequence 
\[
S^{0}(\ker f_{i}) \xrightarrow{} S^{0}E_{i} \xrightarrow{} S^{0}(\coker
f_{i}) \xrightarrow{} S^{1}(\ker f_{i})
\]
in $\stable{\Gamma }$ by Lemma~\ref{lem-stable-cofiber}.  By composing
the maps $S^{0}(\coker f_{i})\xrightarrow{}S^{1}(\ker f_{i})$, we get a
map $S^{0}M\xrightarrow{}S^{k}N$ in $\stable{\Gamma }$.  One can check
that this respects the equivalence relation that defines $\Ext
^{k}_{\Gamma }(M,N)$, and is natural.

Note that this map is an isomorphism when $M=P$ is dualizable, for
then
\[
\stable{\Gamma }(S^{0}P, S^{k}N) \cong \pi _{k}^{P}(S^{0}N) \cong \Ext
^{k}_{\Gamma }(P, N)
\]
by Lemma~\ref{lem-proj-inj}.  Let $\cat{T}$ be the full subcategory
consisting of all $M$ such that 
\[
\alpha _{MN} \mathcolon \Ext ^{k}_{\Gamma }(M,N) \xrightarrow{}
\stable{\Gamma }(S^{0}M, S^{k}N)
\]
is an isomorphism for all $N$ and all $k\geq 0$.  We claim that
$\cat{T}$ is a localizing subcategory.  Indeed, it is clear that
$\cat{T}$ is closed under retracts and coproducts.  To check that
$\cat{T}$ is thick, we note that a short exact sequence 
\[
0 \xrightarrow{} M' \xrightarrow{} M \xrightarrow{} M'' \xrightarrow{} 0
\]
induces a long exact sequence in $\Ext _{\Gamma }^{*}(-,N)$.  Because
short exact sequences are also cofiber sequences in $\stable{\Gamma }$
by~\ref{lem-stable-cofiber}, we also get a long exact sequence in
$\stable{\Gamma }(-, S^{*}N)$.  There is a map between these two long
exact sequences (one must check that $\alpha _{MN}$ is compatible with
the map $\Ext ^{k}_{\Gamma }(M'', N)\xrightarrow{}\Ext ^{k+1}_{\Gamma
}(M',N)$ but the construction of $\alpha _{MN}$ makes this easy to
check).  The five lemma tells that that if two out of $M', M, M''$ are
in $\cat{T}$, so is the third.  
\end{proof}


\begin{thebibliography}{LMSM86}

\bibitem[Ada74]{adams-blue}
J.~F. Adams, \emph{Stable homotopy and generalised homology}, Chicago Lectures
  in Mathematics, University of Chicago Press, Chicago, Ill.--London, 1974,
  x+373 pp.

\bibitem[Bek00]{beke}
Tibor Beke, \emph{Sheafifiable homotopy model categories}, Math. Proc.
  Cambridge Philos. Soc. \textbf{129} (2000), no.~3, 447--475. \MR{1 780 498}

\bibitem[Boa82]{boardman-eight}
J.~M. Boardman, \emph{The eightfold way to {BP}-operations or {$E\sb\ast E$}
  and all that}, Current trends in algebraic topology, Part 1 (London, Ont.,
  1981), CMS Conf. Proc., vol.~2, Amer. Math. Soc., Providence, R.I., 1982,
  pp.~187--226. \MR{84e:55004}

\bibitem[Bor94]{borceux-1}
Francis Borceux, \emph{Handbook of categorical algebra. 1}, Cambridge
  University Press, Cambridge, 1994, Basic category theory. \MR{96g:18001a}

\bibitem[CH02]{hovey-christensen-relative}
J.~Daniel Christensen and Mark Hovey, \emph{Quillen model structures for
  relative homological algebra}, Math. Proc. Cambridge Philos. Soc.
  \textbf{133} (2002), no.~2, 261--293. \MR{1 912 401}

\bibitem[FC90]{faltings-chai}
Gerd Faltings and Ching-Li Chai, \emph{Degeneration of abelian varieties},
  Springer-Verlag, Berlin, 1990, With an appendix by David Mumford.
  \MR{92d:14036}

\bibitem[GH00]{goerss-hopkins-comodules}
Paul~G. Goerss and Michael~J. Hopkins, \emph{Andr\'e-{Q}uillen (co)-homology
  for simplicial algebras over simplicial operads}, Une d\'egustation
  topologique [Topological morsels]: homotopy theory in the Swiss Alps (Arolla,
  1999), Amer. Math. Soc., Providence, RI, 2000, pp.~41--85. \MR{2001m:18012}

\bibitem[Hir02]{hirschhorn}
Philip~S. Hirschhorn, \emph{Model categories and their localizations},
  preprint, available at \texttt{http://www-math.mit.edu/\!$\sim$psh}, August
  4th, 2002.

\bibitem[Hov99]{hovey-model}
Mark Hovey, \emph{Model categories}, American Mathematical Society, Providence,
  RI, 1999. \MR{99h:55031}

\bibitem[Hov01]{hovey-sheaves}
\bysame, \emph{Model category structures on chain complexes of sheaves}, Trans.
  Amer. Math. Soc. \textbf{353} (2001), no.~6, 2441--2457 (electronic). \MR{1
  814 077}

\bibitem[Hov02a]{hovey-barcelona}
\bysame, \emph{Chromatic phenomena in the algebra of ${BP_*BP }$-comodules},
  preprint, 2002.

\bibitem[Hov02b]{hovey-hopf}
\bysame, \emph{Morita theory for {H}opf algebroids and presheaves of
  groupoids}, Amer. J. Math. \textbf{124} (2002), 1289--1318.

\bibitem[HPS97]{hovey-axiomatic}
Mark Hovey, John~H. Palmieri, and Neil~P. Strickland, \emph{Axiomatic stable
  homotopy theory}, Mem. Amer. Math. Soc. \textbf{128} (1997), no.~610, x+114.

\bibitem[HS99a]{hovey-sadofsky-picard}
Mark Hovey and Hal Sadofsky, \emph{Invertible spectra in the ${E}(n)$-local
  stable homotopy category}, J. London Math. Soc. (2) \textbf{60} (1999),
  no.~1, 284--302. \MR{2000h:55017}

\bibitem[HS99b]{hovey-strickland}
Mark Hovey and Neil~P. Strickland, \emph{Morava ${K}$-theories and
  localisation}, Mem. Amer. Math. Soc. \textbf{139} (1999), no.~666, viii+100.

\bibitem[HS02]{hovey-strickland-comodules}
\bysame, \emph{Comodules and {L}andweber exact homology theories}, preprint,
  2002.

\bibitem[JY80]{johnson-yosimura}
David~Copeland Johnson and Zen-ichi Yosimura, \emph{Torsion in
  {B}rown-{P}eterson homology and {H}urewicz homomorphisms}, Osaka J. Math.
  \textbf{17} (1980), no.~1, 117--136. \MR{81b:55010}

\bibitem[Lan76]{land-exact}
P.~S. Landweber, \emph{Homological properties of comodules over
  ${M}{U}_*{M}{U}$ and ${B}{P}_*{B}{P}$}, Amer. J. Math. \textbf{98} (1976),
  591--610.

\bibitem[LMSM86]{lewis-may-steinberger}
L.~G. Lewis, Jr., J.~P. May, M.~Steinberger, and J.~E. McClure,
  \emph{Equivariant stable homotopy theory}, Springer-Verlag, Berlin, 1986,
  With contributions by J. E. McClure. \MR{88e:55002}

\bibitem[May01]{may-triangulated}
J.~P. May, \emph{The additivity of traces in triangulated categories}, Adv.
  Math. \textbf{163} (2001), no.~1, 34--73. \MR{2002k:18019}

\bibitem[MR77]{miller-ravenel}
H.~R. Miller and D.~C. Ravenel, \emph{Morava stabilizer algebras and the
  localization of {N}ovikov's ${E}_{2}$-term}, Duke Math. J. \textbf{44}
  (1977), 433--447.

\bibitem[Pal99]{palmieri-quillen-strat}
John~H. Palmieri, \emph{Quillen stratification for the {S}teenrod algebra},
  Ann. of Math. (2) \textbf{149} (1999), no.~2, 421--449. \MR{1 689 334}

\bibitem[Pal01]{palmieri-book}
John~H. Palmieri, \emph{Stable homotopy over the {S}teenrod algebra}, Mem.
  Amer. Math. Soc. \textbf{151} (2001), no.~716, xiv+172.

\bibitem[Pup79]{puppe-duality}
Dieter Puppe, \emph{Duality in monoidal categories and applications to
  fixed-point theory}, Game theory and related topics (Proc. Sem., Bonn and
  Hagen, 1978), North-Holland, Amsterdam, 1979, pp.~173--185. \MR{82m:55011}

\bibitem[Qui67]{quillen-htpy}
Daniel~G. Quillen, \emph{Homotopical algebra}, Lecture Notes in Mathematics,
  No. 43, Springer-Verlag, Berlin, 1967. \MR{36 \#6480}

\bibitem[Rav86]{ravenel}
D.~C. Ravenel, \emph{Complex cobordism and stable homotopy groups of spheres},
  Academic Press, 1986.

\bibitem[SS00]{schwede-shipley-monoids}
Stefan Schwede and Brooke~E. Shipley, \emph{Algebras and modules in monoidal
  model categories}, Proc. London Math. Soc. (3) \textbf{80} (2000), no.~2,
  491--511. \MR{1 734 325}

\bibitem[Ste75]{stenstrom}
B.~Stenstr\"{o}m, \emph{Rings of quotients}, Die {G}rundlehren der
  mathematischen {W}issenschaften, vol. 217, Springer-{V}erlag, Berlin, 1975.

\end{thebibliography}
\providecommand{\bysame}{\leavevmode\hbox to3em{\hrulefill}\thinspace}
\providecommand{\MR}{\relax\ifhmode\unskip\space\fi MR }
\providecommand{\MRhref}[2]{%
  \href{http://www.ams.org/mathscinet-getitem?mr=#1}{#2}
}
\providecommand{\href}[2]{#2}

\end{document}